\input amstex
\documentstyle{amsppt}
\magnification=1200
\hsize=13.8cm
\catcode`\@=11
\def\NoLogo{\let\logo@\empty}
\catcode`\@=\active
\NoLogo
\def\vgk{\text{\bf VG}_k}
\def\tN{\widetilde N}
\def\tL{\widetilde L}
\def\tM{\widetilde M}
\def\tpi{\widetilde \pi}
\def\hpi{\widehat \pi}
\def\hM{\widehat M}

\def\heat{\lf(\frac{\p}{\p t}-\Delta\ri)}
\def \b {\beta}
\def\i{\sqrt{-1}}
\def\Ric{\text{Ric}}
\def\lf{\left}
\def\ri{\right}
\def\bbar{\bar \beta}
\def\a{\alpha}
\def\ol{\overline}
\def\g{\gamma}
\def\e{\epsilon}
\def\p{\partial}
\def\delbar{{\bar\delta}}
\def\ddbar{\partial\bar\partial}
\def\dbar{\bar\partial}

\def\C{\Bbb C}
\def\R{\Bbb R}
\def\tx{\tilde x}
\def\vp{\varphi}

\def\dbar{\bar\partial}
\def\ba{{\bar\alpha}}
\def\bb{{\bar\beta}}
\def\abb{{\alpha\bar\beta}}

\def\i{\sqrt {-1}}

\def \D {\Delta}
\def\aint{\frac{\ \ }{\ \ }{\hskip -0.4cm}\int}
\documentstyle{amsppt}
\magnification=1200 \hsize=13.8cm \vsize=19 cm

\leftheadtext{Lei Ni  and Luen-fai Tam}
\rightheadtext{Plurisubharmonic functions}
\topmatter
\title{Plurisubharmonic functions and the structure of  complete K\"ahler manifolds with nonnegative curvature}\endtitle

\author{Lei Ni\footnotemark and  Luen-Fai Tam\footnotemark}\endauthor
\footnotetext"$^{1}$"{Research partially supported by NSF grant
DMS-0203023, USA.} \footnotetext"$^2$"{Research partially
supported by Earmarked Grant of Hong Kong \#CUHK4032/02P.}
\address
Department of Mathematics, University of California, San Diego, La Jolla,
CA 92093
\endaddress
\email{
lni\@math.ucsd.edu}
\endemail

\address
Department of Mathematics, The Chinese University of Hong Kong,
Shatin, Hong Kong, China
\endaddress
\email{lftam\@math.cuhk.edu.hk}
\endemail

\affil
{
University of California, San Diego\\
The Chinese University of Hong Kong
}
\endaffil

\date  February  2003, revised in March 2003\enddate

\abstract
In this paper, we study  global properties of
continuous plurisubharmonic functions on complete noncompact
K\"ahler manifolds with nonnegative bisectional curvature and
their applications to the structure of such manifolds. We prove that
continuous plurisubharmonic functions with    reasonable  growth rate  on
such manifolds can be approximated by smooth plurisubharmonic functions
through the heat flow deformation. Optimal Liouville type theorem for the
plurisubharmonic functions as well as a splitting theorem in terms of
harmonic functions and  holomorphic functions are established. The results are
then applied to prove several structure theorems on complete noncompact
K\"ahler manifolds with nonnegative bisectional or sectional curvature.
\endabstract

\endtopmatter
\document
\subheading{\S0 Introduction}\vskip .2cm

In this paper,  we are interested in the class of  complete
noncompact K\"ahler manifolds with nonnegative holomorphic
bisectional curvature. We shall first give a detailed study on the
properties of heat flow with plurisubharmonic functions as initial
data. Then we shall use the results   to prove a Liouville theorem on plurisubharmonic
functions and a splitting theorem related to harmonic and
holomorphic functions. All these results will then be applied to
obtain structure theorems on K\"ahler manifolds with nonnegative
sectional or holomorphic bisectional curvature.

One  motivation of the present work is a  program proposed by Yau
\cite{Y 3, p. 622} on the study of parabolic manifolds:  ``The
question is to demonstrate that every noncompact K\"ahler manifold
with positive bisectional curvature is biholomorphic to the
complex euclidean space. If we only assume that the nonnegativity
of the bisectional curvature, the manifold should be biholomorphic
to a complex vector bundle over a compact Hermitian symmetric
space.'' As pointed out in \cite{Y 3}, an important reason for  this program comes from the celebrated results of Cheeger-Gromoll
\cite{CG 2} and Gromoll-Meyer \cite{GM} on complete noncompact
Riemannian manifolds with nonnegative or positive sectional
curvature. It is also motivated by the work of Greene-Wu \cite{GW
2} on the Steinness of K\"ahler manifolds.   In both cases, a key
ingredient is to study Busemann functions.

In \cite{CG 2} it was proved  that the Busemann function (with
respect to all geodesic rays from a fixed point) on a complete
noncompact Riemannian manifold with nonnegative sectional
curvature    is Lipschitz continuous, convex and is an exhaustion
function. Then it was proved that a complete noncompact Riemannian
manifold with nonnegative sectional curvature is diffeomorphic to
the normal bundle over a compact totally geodesic submanifold
without boundary, which is totally convex and is called the `soul' of
the manifold.

On a  K\"ahler manifold with nonnegative holomorphic bisectional
curvature, even though the Busemann function is no longer  convex,
it is still plurisubharmonic. This was proved by Wu \cite{W 1}.
Moreover, it is in fact strictly plurisubharmonic at the point
where the holomorphic bisectional curvature is positive. Using
this fact, it was proved  by Greene-Wu  that if  the manifold has
nonnegative sectional curvature and positive holomorphic
bisectional curvature, then  it is Stein because in this case the
Busemann function is also an exhaustion function, see \cite{GW
1--4, W 1-2} for more results. In the proof, it was first shown
that a continuous {\it strictly} plurisubharmonic function can be
approximated uniformly by a smooth one. The result of  Grauert
\cite{G} can then be applied to conclude that the manifold is
Stein.

If the manifold has nonnegative holomorphic bisectional curvature,
then the Busemann function is only (continuous) plurisubharmonic instead of
strictly plurisubharmonic. In order to use the Busemann function,
it is desirable to approximate it by a smooth one.  In general, it
seems unlikely that a continuous plurisubharmonic function can be
approximated by $C^\infty$-plurisubharmonic functions. However, we
shall prove the following rather general results on the solution
of the heat equation with continuous plurisubharmonic function as
initial data (see Theorem 3.1).

\proclaim{Theorem 0.1} Let $M^m$ be a complete noncompact K\"ahler
manifold with nonnegative holomorphic bisectional curvature and
let  $u$ be a continuous plurisubharmonic function on $M$
satisfying
$$
|u|(x)\le C \exp(ar^2(x)) \tag0.1
$$
for some positive constants
$a$, $C$, where $r(x)$ is the
distance of $x$ from a fixed point. Let $v$ be the solution of the
heat equation with initial data $u$. There exists $T_0>0$
depending only on $a$ and there exists $T_0>T_1>0$ such that the
following are true.
\roster
\item"{(i)}" For $0<t\le T_0$,
$v(\cdot,t)$ is a smooth plurisubharmonic function.
\item"{(ii)}" Let
$$
\Cal K(x,t)=\{w\in T^{1,0}_x(M)|\  v_\abb(x,t) w^\a=0,\text{\rm \ for all $\b$}\}
$$
be the null space of $v_\abb(x,t)$. Then for any $0<t<T_1$, $\Cal
K(x,t)$ is   distribution on $M$. Moreover the distribution is
invariant under parallel translations.
\item"{(iii)}" If the
holomorphic bisectional curvature is positive at some point, then
$v(x,t)$ is strictly plurisubharmonic for all $0<t<T_1$.
\endroster
\endproclaim
In particular,  if $u$ is a continuous plurisubharmonic function
satisfying (0.1), then it can be approximated by smooth
plurisubharmonic functions uniformly on compact subsets. In
application, we shall make use of the properties of $v(x,t)$ in
the above theorem rather than the result on approximation.

The proof of Theorem 0.1  relies on a  general maximum principle
for Hermitian symmetric (1,1) tensor $\eta$ which  satisfies a linear
heat equation  on a complete noncompact K\"ahler manifold with
nonnegative holomorphic bisectional curvature. We obtain a maximum
principle for $\eta$ under some weak growth conditions on the rate
of the average of $||\eta||$, the norm  of $\eta$, over geodesic
balls. Since there is no pointwise bound on $||\eta||$, we shall
apply an indirect cutoff argument together with careful estimates
on the solutions of the heat equation through extensive uses of
the fundamental work of Li and Yau \cite{LY}.

In \cite{N 1},
the first author raised  the following question:

\medskip

{\it On a complete noncompact K\"ahler manifold with nonnegative Ricci curvature,
is  a plurisubharmonic function of sub-logarithmic growth  a constant?}

\medskip

\noindent It is well-known that for the complex Euclidean space $\C^m$, the
answer is positive.  An affirmative answer to the above question is also a natural
analogue, for plurisubharmonic functions, of Yau's Liouville theorem [Y 1] for
positive harmonic functions on Riemannian manifolds with nonnegative Ricci
curvature.
An immediate application of the Theorem 0.1 is to  give an affirmative answer to the above question on K\"ahler manifolds with nonnegative holomorphic bisectional curvature. Namely, we have   the following (see Theorem 3.2):

\proclaim{Theorem 0.2} Let $M$ be a complete noncompact K\"ahler manifold
with nonnegative holomorphic bisectional curvature. Let $u$ be a continuous
plurisubharmonic function on $M$. Suppose that
$$
\limsup_{x\to \infty} \frac{u(x)}{\log r(x)}=0.
$$
Then $u$ must be a constant.
\endproclaim

Using Theorem 0.2 we obtain the following interesting results (see Theorem 4.1):\proclaim{Theorem 0.3} Let $M^m$ be a complete noncompact K\"ahler
manifold with nonnegative holomorphic bisectional curvature.
Suppose $f$ is a nonconstant harmonic function on $M$ such that
$$
\limsup_{x\to\infty}\frac{|f(x)|}{r^{1+\e}(x)}=0, \tag0.2
$$
for any $\e>0$, where $r(x)$ is the distance of $x$ from a fixed
point. Then  $f$ must be of linear growth and  $M$ splits
  isometrically as $\widetilde{M}\times\R$.
 Moreover the universal cover $\overline {M}$ of $M$  splits
  isometrically and holomorphically as $\widetilde{M'}\times\C$,
where $\widetilde M'$ is a complete K\"ahler
  manifold with nonnegative holomorphic bisectional curvature. Suppose that
there exists a nonconstant holomorphic function $f$ on $M$ satisfying (0.2). Then $M$ itself splits as $\widetilde{M}\times\C$.
\endproclaim
 A well-known result in \cite{Y 1, CY} says that if the growth rate of a harmonic function on a complete noncompact Riemannian manifold with nonnegative Ricci curvature is `close' to constant functions, namely if it is of sublinear growth, then it must be constant. Similar to this, the first result of Theorem 0.3 says that on a complete noncompact K\"ahler manifold with nonnegative holomorphic bisectional curvature, if the growth rate of a harmonic function is `close' to linear, then it must be of linear growth. On the other hand, for any $\delta>0$,
the `round off' cones with metrics
$dr^2+r^2\ ds^2_{S^1(1+\delta)}$, where $S^1(\frac{1}{\sqrt{1+\delta}})$
is the circle with radius $\frac{1}{\sqrt{1+\delta}}$, supports harmonic functions of growth $r^{1+\delta}(x)$.

    One might also want to compare the splitting result in Theorem 0.3 with some  previous related results  in \cite{CG 1, L 2, CCM}. In \cite{GG 1}, it was proved that if a complete noncompact Riemannian manifold with nonnegative Ricci curvature contains a line then a factor $\R$ can be splitted isometrically. In \cite{L2}, it was proved that if a complete noncompact K\"ahler manifold with nonnegative Ricci curvature with complex dimension $m=n/2$ supports $n+1$ independent linear growth harmonic functions, then it is isometric and holomorphic to $\C^m$. In \cite{CCM}, Li's result was generalized to the Riemannian case, and the conclusion is that the manifold is isometric to Euclidean space. In \cite{CCM}, result  on the splitting of the tangent cone in terms of linear growth harmonic functions was obtained.

The second part
 of the paper is to  study the structure of  complete noncompact  K\"ahler manifolds with nonnegative holomorphic bisectional curvature.
The main tool is to use the heat flow with the Busemann functions as initial data.  As mentioned above, the Busemann function is a continuous plurisubharmonic function on a complete noncompact K\"ahler manifold with nonnegative holomorphic bisectional curvature. Hence Theorem 0.1 will be very useful. It turns out that Theorem 0.3 will  be useful in the study too.

  Before we state our next result, let
us first introduce some conditions on a K\"ahler manifold. The first one
is on  the  growth rate of   volumes of geodesic balls.   $M$ is said to satisfy
($\vgk$) for $k>0$, if there exists a constant $C>0$ such that
$$
V_o(r)\ge Cr^{k}\tag"($\vgk$)"
$$
for all $r\ge 1$. The other two conditions are on the decay of the
curvature. Suppose $M$ has nonnegative scalar curvature $\Cal R$.
$M$ is said to satisfy the curvature decay condition ({\bf CD}) if
there exists a constant $C>0$ (which might depends on $o$) such
that
$$
\aint_{B_o(r)}\Cal R\le \frac{C}{r}\tag"({\bf CD})"
$$
for all $r>0$.  $M$ is said to satisfy the fast curvature decay
condition ({\bf FCD}) if there is a constant $C>0$, so that
$$
\int_0^rs\lf(\aint_{B_o(s)}\Cal R(x)dx\ri)ds\le C\log(r+2)\tag"({\bf FCD})"
$$
for all $r>0$. ({\bf FCD})  means that the average of the scalar
curvature decays quadratically in the integral sense. Hence it is
stronger than ({\bf CD}). Our next result is the following
splitting theorem (see Theorem 4.2):

\proclaim{Theorem 0.4} Let $M^m$ be a complete noncompact K\"ahler
manifold with nonnegative holomorphic bisectional curvature.

\roster
\item"{(i)}" Suppose $M$ is simply connected, then $M=N\times M'$
holomorphically and isometrically, where $N$ is a compact simply connected K\"ahler manifold, $M'$ is a complete noncompact K\"ahler
manifold and both $N$ and $M'$ have  nonnegative holomorphic bisectional curvature.
Moreover, $M'$ supports a smooth strictly plurisubharmonic
function with bounded gradient and   satisfies {\rm($\vgk$)} and {\rm({\bf
CD})}, where $k$ is the complex dimension of $M'$. If, in addition,
$M$ has nonnegative sectional curvature outside a compact set,
then $M'$ is also Stein.

\item"{(ii)}" If the holomorphic bisectional  curvature of $M$ is
positive at some point, then $M$ itself supports a smooth strictly
plurisubharmonic function with bounded gradient and satisfies
{\rm($\text{\bf VG}_m $)} and {\rm({\bf CD})}, where $m$ is the complex
dimension of $M$. If, in addition, $M$ has nonnegative sectional
curvature outside a compact set, then $M$ is also Stein.
\endroster
\endproclaim

The conclusion on the volume growth in the first statement in (ii)
was first proved in   \cite{CZ 2} and the conclusion on curvature
decay  is a generalization of a result in \cite{CZ 2}.  The last
statement in (ii) is a generalization of a result in \cite{W 1}.
Note that by \cite{M 2, HSW} (see also \cite{CC}), $N$ in the
theorem is a compact Hermitian symmetric manifold, but we shall
not use this fact in the proof. Note also that $N$ may not be
present.

An immediate consequence of Theorem 0.4 is on the Steinness of complete noncompact K\"ahler manifolds with nonnegative holomorphic bisectional curvature. Recall that a complete noncompact Riemannian manifold of dimension $n$ with nonnegative Ricci curvature is said to have maximum volume growth if $V_x(r)\ge Cr^n$ for some positive constant $C$ for all $x$ and $r$.  A result in \cite{Sh}  states that the Busemann
function on a complete noncompact manifold with nonnegative Ricci
curvature and with maximum volume growth is an exhaustion
function. Using this, we prove as a corollary to Theorem 0.4  that a complete noncompact K\"ahler
manifold with nonnegative holomorphic bisectional curvature and
maximum volume growth is Stein. Here we assume neither that the
holomorphic bisectional curvature is positive, which implies that
the Busemann function is strictly plurisubharmonic,  nor any
curvature decay conditions  as in \cite{CZ 1}. We also prove the
Steinness  for the case that the manifold has a pole.  This answers a question raised
in \cite{W 2, page 255} affirmatively. Recall that a Riemannian manifold is said to have a pole if there is a point $p$ in the manifold such that the exponential map at $p$ is a diffeomorphism.

To study   $M'$ (or $M$ in Theorem 0.4(ii)) further, we obtain the
following (see Theorem 4.3):

 \proclaim{Theorem 0.5} Let $M^m$ be a complete noncompact K\"ahler
manifold with nonnegative holomorphic bisectional curvature.
Assume that $M$ supports a smooth strictly plurisubharmonic
function $u$ on $M$ with bounded gradient.

\roster
\item"{(i)}" If $M$ is simply connected, then $M=\C^\ell\times M_1\times
M_2 $ isometrically and holomorphically for some $\ell\ge0$, where
$M_1$ and $M_2$ are complete noncompact K\"ahler manifold with
nonnegative holomorphic bisectional curvature such that any
polynomial growth holomorphic function on $M$ is independent of
the factor $M_2$, and any linear growth holomorphic function is
independent of the factors $M_1$ and $M_2$. Moreover, $M_1$
supports a strictly plurisubharmonic function of logarithmic growth
and  satisfies {\rm({\bf FCD})} and {\rm($\text{\bf VG}_{a}$)}, for any
$a<k+1$,  where $k=\dim_\C M_1$.

\item"{(ii)}" Suppose the holomorphic bisectional curvature of $M$ is
positive at some point, then either $M$ has no nonconstant
polynomial growth holomorphic function or $M$ itself satisfies
{\rm({\bf FCD})} and {\rm($\text{\bf VG}_{a}$)}, for any $a<m+1$.
\endroster
\endproclaim
There is an open question on whether the ring of polynomial growth
holomorphic functions on a complete noncompact K\"ahler manifold
with nonnegative curvature is finitely generated, see \cite{Y 4,
p. 391}. By Theorems 0.4 and 0.5, in order to study polynomial
growth holomorphic functions on  a manifold with nonnegative
holomorphic bisectional curvature    which is either simply
connected or has positive holomorphic bisectional curvature at
some point, we may assume that $M$ satisfies  the fast curvature
decay condition ({\bf FCD}) and the volume growth condition
($\text{VG}_a$) for any $a<m+1$.

Together with the   $L^2$ estimates \cite{H\"o} and the mean value
inequality \cite{LS} Theorems 0.4 and 0.5 also imply that a simply
connected complete noncompact K\"ahler manifold $M$ with
nonnegative holomorphic bisectional curvature supports many
nontrivial holomorphic functions. Namely,   $M$ is a product of a
compact Hermitian symmetric manifold, a complex Euclidean space, a
complete manifold $M_2$  and a complete manifold $M_1$ such that
each point of $M_2$ has   local coordinate functions which are the
restriction of global holomorphic functions with exponential
growth of order $\le 1$ in the sense of Hadamard, and each point
$M_1$ has   local coordinate functions which are the restriction
of global holomorphic functions with polynomial growth.

The results in the theorems on the decay rate of the average of
the scalar curvature are  related to the work of Shi \cite{S} on
the long time existence of the K\"ahler-Ricci flow, see also
\cite{Y 3}. Theorems 0.4 and 0.5 also imply some  uniformization
type results when the volume growth of the manifold is small, see
Corollary 4.3. Namely a simply-connected complete K\"ahler
manifold with nonnegative bisectional curvature and slow volume
growth is biholomorphic to the product of the complex line with a
compact Hermitian- symmetric spaces.

Next we shall study K\"ahler manifolds  with nonnegative
holomorphic bisectional curvature whose Busemann functions are   exhaustion functions. We also assume that the universal cover of the manifold
does not contain de Rham Euclidean factors. This class of manifolds
contains the manifolds which have nonnegative sectional curvature
outside a compact set and positive Ricci curvature somewhere.
Without assuming that the manifold is simply connected, one can
describe the structure of $M$ in a rather explicit way, see Theorems 5.1.

\proclaim{Theorem 0.6}   Let $M$ be a complete noncompact K\"ahler
manifold   with nonnegative holomorphic bisectional curvature such
that the Busemann function is an exhaustion function. Suppose the
universal cover $\tM$ has no Euclidean factor. Then $\tM=\tN\times \tL$ where $\tN$ is a compact Hermitian symmetric manifold and $\tL$ is Stein. Moreover, $M$ is a holomorphic and Riemannian fibre bundle with fibre $\tN$ over a Stein manifold $\hM$ with nonnegative holomorphic bisectional curvature such that $\tM$ is covered by $\tL$.
\endproclaim

Using this  structure result, Fangyang Zheng \cite{Z 2} proves that  if in addition that $M$  has nonnegative sectional curvature everywhere, $M$ is in fact simply-connected and $M=N\times L$, where $N$ is compact  $L$ is a Stein manifold and is diffeomorphic to $R^{2l}$ where $l=\dim_{\C} L$. From this, one can prove that a complete noncompact K\"ahler manifold with nonnegative sectional curvature is a holomorphic and Riemannian fibre bundle over $\C^k/\Gamma$ for some discrete subgroup of the holomorphic isometry group of $\C^k$, with fibre $N\times L$ with the structures as above. The authors are grateful to   Fangyang Zheng for  allowing us to include his results and proofs in this work, see Theorem 5.2 and Corollary 5.1.

The results are motivated by the work of  Takayama \cite{Ta},
where he proved that if $M$ is a complete noncompact K\"ahler
manifold with nonnegative holomorphic bisectional curvature and
negative canonical line bundle and if $M$ supports a continuous
plurisubharmonic exhaustion function, then $M$ has a structure of
holomorphic fibre bundle over a Stein manifold whose fibre is
biholomorphic to some compact Hermitian symmetric manifold.
Obviously, our assumptions are stronger. However, in Theorem 0.6,
the structure of the manifold is described more explicitly.
Moreover, our proof is rather elementary and does not appeal to
the result of \cite{M 2} for example.

Finally, the methods of our study on the heat equation and the
maximum principle can be applied to obtain the following result
which is related Theorem 0.4 and in particular to the condition
({\bf FCD}).

\proclaim{Theorem  0.7} Let $M$ be a complete noncompact K\"ahler
manifold with nonnegative holomorphic bisectional curvature. Then
$M$ is flat if
$$
\int_0^rs\lf(\aint_{B_o(s)}\Cal R(y)\, dy\ri)\, ds=o(\log r)
$$
provided   that

$$
\liminf_{r\to\infty}\lf[\exp\lf(-ar^2\ri)\int_{B_o(r)}\Cal R^2\ri]<\infty$$
for some $a>0$.
where $\Cal R$ is the    scalar curvature  of $M$.
\endproclaim
For previous results in this
 direction, see \cite{MSY, N 1, NST 1, CZ 1}.
  One of the main ideas is to solve the Poincar\'e-Lelong equation under rather weak conditions and then apply Theorem 0.2.
  The solution of the Poincar\'e-Lelong equation may has independent interest. See previous works \cite{MSY, NST 1} on this problem.

Recently, Wu and Zheng   \cite {WZ  1-2} prove some    splitting results  on    K\"ahler manifolds   with  nonnegative  or with non-positive  holomorphic bisectional curvature in terms of the rank of the Ricci tensor. In their works,  the metric is assumed to be real analytic.

We organize the paper as follows: in \S1 we study the solution of
the heat equation; in \S2 a maximum principle for Hermitian
symmetric (1,1) tensor is given; we then apply the results to
study the solution of the heat flow with continuous
plurisubharmonic initial data in \S3, a Liouville theorem for
plurisubharmonic functions is also proved there; in \S4--\S6, we
shall discuss the structure of K\"ahler manifolds with nonnegative
holomorphic bisectional curvature; a solution to the
Poincar\'e-Lelong equation will also be given in \S6.

The authors would like to thank Professors  Laszlo Lempert, Hing Sun Luk, Shigeharu Takayama, Hung-Hsi Wu and Fangyang Zheng for some useful discussions. They also would like to thank Professors Peter Li and Richard Schoen for their interest in this work.

\input amstex
\documentstyle{amsppt}
\magnification=1200 \hsize=13.8cm \catcode`\@=11
\def\NoLogo{\let\logo@\empty}
\catcode`\@=\active \NoLogo
\def\heatt{\lf (\Delta-\frac{\p}{\p t}\ri)}

\def\heat{\lf(\frac{\p}{\p t}-\Delta\ri)}
\def\fp{f_+}
\def\fm{f_-}
\def \b {\beta}
\def\i{\sqrt{-1}}
\def\Ric{\text{Ric}}
\def\lf{\left}
\def\ri{\right}
\def\bbar{\bar \beta}
\def\a{\alpha}
\def\ol{\overline}
\def\g{\gamma}
\def\e{\epsilon}
\def\p{\partial}
\def\delbar{{\bar\delta}}
\def\ddbar{\partial\bar\partial}
\def\dbar{\bar\partial}

\def\C{\Bbb C}
\def\R{\Bbb R}
\def\tx{\tilde x}
\def\vp{\varphi}

\def\tN{\tilde\nabla}

\def\dbar{\bar\partial}
\def\ba{{\bar\alpha}}
\def\bb{{\bar\beta}}
\def\abb{{\alpha\bar\beta}}

\def\i{\sqrt {-1}}

\def \D {\Delta}
\def\aint{\frac{\ \ }{\ \ }{\hskip -0.4cm}\int}
\documentstyle{amsppt}
\vsize=19.0 cm

\subheading{\S1 Preliminary results}\vskip .2cm

In this section, we shall derive some basic results on the
solutions to the heat equation on a complete noncompact manifold
with nonnegative Ricci curvature. These results will be used in
later sections regularly. Specifically, we shall show that the
Cauchy problem (1.6), which shall be defined  in the following,
can be solved, where only the average growth rate of the initial
data over geodesic balls is assumed. This condition is useful in
applications because in many cases a continuous plurisubharmonic
function can only be approximated by a smooth function without
point-wise estimations on  the norms of the complex Hessians.
However, Lemma 1.6 below shows that they  can be estimated  in the
average sense. Corollary 1.4 and Lemma 1.4 will be used to keep
track of the behaviors of the functions which approximate  the
Buesmann function through the heat flow. These are important in
the study of the structures of the manifolds.

We always assume that  $M^n$ is a complete noncompact Riemannian
manifold with nonnegative Ricci curvature within this section. Let
$H(x,y,t)$ be the heat kernel of $M$ and let $o\in M$ be a fixed
point. Denote the average of a function $f$ over $B_x(r)$ by
$\aint_{B_x(r)}f$. In this work, we shall make extensive uses of
the fundamental work on the heat kernel estimates of Li and Yau
\cite{LY}. We start with a $L^p$-estimate on the nonnegative
 solution to the heat equation.

\proclaim{Lemma 1.1} Let $f\ge 0$ be a function on a complete noncompact Riemannian manifold $M^n$ with nonnegative Ricci curvature and let
$$
u(x,t)=\int_MH(x,y,t)f(y)dy.
$$
Assume that $u$ is defined on $M\times[0,T]$ for some $T>0$ and that for $0<t\le T$,
$$
\lim_{r\to\infty}\exp\lf(-\frac{r^2}{20t}\ri)\int_{B_o(r)}
f=0.\tag1.1
$$
Then for any $r^2\ge t>0$, and $p\ge 1$,
$$
\aint_{B_o(r)}u^p(x,t)dx\le C(n,p)\lf[\aint_{B_o(4r)}f^p(x)dx+t^{-p}\lf(\int_{4r}^\infty \exp\lf(-\frac{s^2}{20 t}\ri)s\aint_{B_o(s)}fds\ri)^p\ri].
$$
\endproclaim
\demo{Proof} For any $p\ge 1$ and $r\ge \sqrt t$,
$$
\split
\int_{B_o(r)} u^p(x,t)dx &=\int_{B_o(r)}\lf(\int_MH(x,y,t)f(y)dy\ri)^p dx\\
&\le C(p)\bigg[\int_{B_o(r)}\lf(\int_{B_o(4r)}H(x,y,t)f(y)dy\ri)^pdx\\
&\qquad+\int_{B_o(r)}\lf(\int_{M\setminus B_o(4r)}H(x,y,t)f(y)dy\ri)^pdx.\bigg]
\endsplit\tag1.2
$$
Now for $x\in B_o(r)$ and $y\notin B_o(4r)$, we have $r(x,y)\ge
3/4r(y)$.  By the estimates of the heat kernel of Li and Yau
\cite{LY, p.176},  we have:
$$
\split
\int_{M\setminus B_o(4r)}&H(x,y,t)f(y)dy\\
&\le C_1 \int_{M\setminus B_o(4r)}\frac1{V_x(\sqrt t)}\exp(-\frac{r^2(x,y)}{5t})f(y)dy\\
&\le C_2\int_{M\setminus B_o(4r)}\frac1{V_x(r+\sqrt t)}\cdot\lf(\frac{r+\sqrt t}{\sqrt t}\ri)^n \exp(-\frac{r^2(x,y)}{5t})f(y)dy\\
&\le\frac{C_2}{V_o(\sqrt t)}\cdot\lf(\frac{r+\sqrt t}{\sqrt t}\ri)^n\int_{M\setminus B_o(4r)}\exp(-\frac{r^2(x,y)}{5t})f(y)dy\\
& \le\frac{C_2}{V_o(\sqrt t)}\cdot\lf(\frac{r+\sqrt t}{\sqrt t}\ri)^n\int_{4r}^\infty \exp(-\frac{s^2}{10t})\lf(\int_{\p B_o(s)}f\ri) ds\\
&\le \frac{C_2}{10V_o(\sqrt t)}\cdot\lf(\frac{r+\sqrt t}{\sqrt t}\ri)^n \int_{4r}^\infty \exp(-\frac{s^2}{10t})\lf(\int_{  B_o(s)}f\ri) d\lf(\frac{s^2}t\ri)\\
& \le C_3\lf(\frac{r+\sqrt t}{\sqrt t}\ri)^n \int_{4r}^\infty \frac{V_o(s)}{V_o(\sqrt t)}\cdot\exp(-\frac{s^2}{10t})\lf(\aint_{  B_o(s)}f\ri) d\lf(\frac{s^2}t\ri)\\
&\le C_4t^{-1}\lf[\int_{4r}^\infty\lf(\frac{s}{\sqrt t}\ri)^{2n}\exp(-\frac{s^2}{10t})s\aint_{  B_o(s)}fds\ri]\\
&\le C_5t^{-1}  \lf[\int_{4r}^\infty \exp(-\frac{s^2}{20t})s\aint_{  B_o(s)}fds\ri]
\endsplit\tag1.3
$$
for some constants $C_1-C_5$ depending only on $n$. Here we have
used the volume comparison and the assumption (1.1) when we
perform integration by parts in the fifth inequality.

On the other hand,  by H\"older inequality and   the fact that $\int_MH(x,y,t)dy=1$, we have
$$
\lf(\int_{B_o(4r)}H(x,y,t)f(y)dy\ri)^p\le\int_{B_o(4r)}H(x,y,t)f^p(y)dy.
$$
Hence
$$
\split
\int_{B_o(r)}\lf(\int_{B_o(4r)}H(x,y,t)f(y)dy\ri)^pdx&\le \int_{B_o(r)}\int_{B_o(4r)}H(x,y,t)f^p(y)dy\,dx\\
&\le \int_{B_o(4r)}f^p(y)\lf(\int_{B_o(r)}H(x,y,t)dx\ri) dy\\
&\le \int_{B_o(4r)}f^p(y)dy.
\endsplit\tag1.4
$$
The lemma follows from (1.2)--(1.4).
\enddemo

 Let $u$ be a continuous function on $M$ such that
$$
\aint_{B_o(r)}|u|(x)\, dx\le \exp(ar^2+b)\tag1.5
$$
for some positive constant $a>0$ and $b>0$. Consider the following initial value problem
$$
\cases & \heatt v(x,t)=0\\
&\quad v(x,0)=u(x).
\endcases\tag1.6
$$
\proclaim{Lemma 1.2} The initial value problem (1.6) has a solution on $M\times[0,\frac1{40a}]$. Moreover, for $(x,t)\in M\times(0,\frac1{40a}]$,
$$
v(x,t)=\int_MH(x,y,t)u(y)dy,
$$
where $H(x,y,t)$ is the heat kernel of $M$
\endproclaim
\demo{Proof}  For $j\ge1$,  let $0\le \varphi_j\le1$ be a smooth cutoff function such that $\varphi_j\equiv 1$ on $B_o(j)$ and $\varphi_j\equiv0$ on $B_o(2j)$. Let $u_j=\varphi_j u$. Then $u_j$ is continuous with compact support. Hence one can solve (1.6) with initial value $u_j$ for all time. The solution $v_j$ is given by
$$
v_j(x,t)=\int_M H(x,y,t)u_j(y)dy\tag1.7
$$
for $(x,t)\in M\times(0,\infty)$. By Lemma 1.1, for any   $0<t\le \frac1{40a}$ and for any $r\ge\sqrt t$,
$$
\split
\aint_{B_o(r)}|v_j|(x,t)&\le \aint_{x\in B_o(r)}\lf(\int_M H(x,y,t)|u_j|(y)dy\ri)  dx\\
&\le C_1 \lf[\aint_{B_o(4r)}|u_j|+t^{-1}\int_{4r}^\infty \exp(-\frac{s^2}{20 t})s\aint_{B_o(s)}|u_j|ds\ri]\\
&\le C_2\lf[\aint_{B_o(4r)}|u_j|+e^bt^{-1}\int_{4r}^\infty \exp(-\frac{s^2}{20 t}+as^2)s ds\ri]\\
&\le C_3e^b\lf[\exp(16ar^2)+ \int_{4r}^\infty \exp(-\frac{s^2}{40 t})  d\lf(\frac{s^2}{t}\ri)\ri]\\
&\le C_4e^b\lf(\exp(16ar^2)+1\ri)
\endsplit\tag1.8
$$
where $C_1-C_4$ are constants depending only on $n$. Since $|v_j|$
are subsolutions of  the heat equation and $|u_j|\le |u|$, by
\cite{LT, Theorem 1.2} and (1.7) for $R^2\ge 1/(40a)$, we have
that
$$
\sup_{B_o(\frac 12 R)\times[0,\frac1{40a}]}|v_j|\le C_5\lf[\exp(16 aR^2+b) +\sup_{B_o(R)}|u|\ri] \tag1.9
$$
for some constant $C_5$ depending only on $n$.
From this, it is easy to see that after passing to a subsequence, $v_j$ together their derivatives converge  uniformly on compact sets on $M\times(0,\frac1{40a}]$ to a solution $v$ of the heat equation. Moreover, for any $(x,t)\in M\times(0,\frac1{40a}]$, as in (1.3) we have
$$
\split
\lf|\int_MH(x,y,t)u(y)dy-v_j(x,t)\ri|&=
\lf|\int_MH(x,y,t)\lf(u (y)-u_j(y)\ri)dy\ri|\\
&\le \int_{M\setminus B_o(j)}H(x,y,t)|u|(y)dy\\
&\le C_6
\int_{j}^\infty \exp(-\frac{s^2}{20t})s\aint_{\p B_o(s)}|u|ds\\
&\le C_6\int_j^\infty \exp (-\frac{s^2}{40t})d\lf(\frac{s^2}t\ri)\\
&\le C_6\int_{\frac{j^2}t}^\infty\exp(-\frac1{40}\tau)d\tau
\endsplit
$$
for some positive constants $C_6$. Here we have used the Harnack
inequality \cite{LY, P. 168},  the assumption (1.5) on $u$ and the
fact that $t\le \frac{1}{40 a}.$ Hence it is easy to see that
$$
v(x,t)=\int_M H(x,y,t)u(y)dy
$$
and $v(x,0)=u(x)$.
 \enddemo

In the next lemma, we shall obtain an estimate of the growth rate of $v(x,t)$ for fixed $t$ in terms of the growth rate of $u$.

\proclaim{Lemma 1.3} Let $u$ and $v$ be as in Lemma 1.2. Then for any $1>\e>0$, there exists a constant  $C=C(n,\e,a,b)$ depending only on $n$, $\e$, $a$ and $b$,  and there exists $\frac{1}{40a}>T_0>0$ depending only on $a$ and $\e$, such that for all $x\in M\times(0,T_0]$ with $r^2(x)\ge  T_0 $,
$$
\lf|v(x,t)-\int_{B_x(\e r)}H(x,y,t)u(y)dy\ri|\le C(n,\e,a,b)
$$
where $r=r(x)$.
\endproclaim
\demo{Proof} Let $x\in M$ and let $r=r(x)$. It is easy to see that
for $s\ge \e r$
$$
\aint_{B_x(s)}|u|(y)\, dy \le
C_1\aint_{B_o\lf((1+\e^{-1})s\ri)}|u|(y)\, dy
$$
for some constant $C_1$ depending only on $n$ and $\e$. Hence if
$T_0>0$ is small enough, depending only on $\e$ and $a$, then for
$0<t\le T_0$,  as before by \cite{LY, p. 176} we have
$$
\split
\int_{M\setminus B_x(\e r)}H(x,y,t)|u|(y)dy&\le \frac{C_2}{V_x(\sqrt t)}
\int_{\e r}^\infty \exp\lf(-\frac{s^2}{5t}\ri) \lf(\int_{B_x(s)}|u|(y)\, dy\ri) d\lf(\frac{s^2}{t}\ri)\\
&\le C_3\int_{\e r}^\infty \lf(\frac{s}{\sqrt t}\ri)^n\exp\lf(-\frac{s^2}{5t}\ri)
\lf(\aint_{B_o(\lf(1+\e^{-1})s\ri)}|u|\, dy\ri) d\lf(\frac{s^2}{t}\ri)\\
&\le C_3\int_{\e r}^\infty \lf(\frac{s}{\sqrt t}\ri)^n\exp\lf[-\frac{s^2}{5t}+a\lf(1+\e^{-1}\ri)^2s^2\ri] d\lf(\frac{s^2}{t}\ri)\\
&\le C_4
\endsplit
$$
for some constants $C_2$, $C_3$ $C_4$ depending only on $n$, $\e$,
$a$ and $b$. From this the lemma follows.
\enddemo

\proclaim{Corollary 1.4} With the same assumptions and notations as in Lemma 1.3, let $ C(n,\e,a,b)$ be the constant in the lemma. Then for $x\in M$ with $r=r(x)\ge \sqrt{T_0}$ such that $u\ge 0$ on $B_x(\e r)$, then for any $0\le t<T_0$
$$
-C(n,\e,a,b)+C_1\inf_{B_x(\e r)}u\le v(x,t)\le C(n,\e,a,b)+\sup_{B_x(\e r)}u
$$
for some positive constant $C_1$ depending only on $n$ and $\e$.
\endproclaim
\demo{Proof} By Lemma 1.3, since $\int_MH(x,y,t)\,dy=1$, we have
$$
v(x,t)\le C(n,\e,a,b)+\int_{B_x(\e r)}H(x,y,t)u(y)\,dy\le C(n,\e,a,b)+\sup_{B_x(\e r)}u.
$$
 On the other hand, by
the lower bound estimate of the heat kernel of Li-Yau \cite{LY, p.182} and Lemma 1.3,   we have
that
$$
\split
 v(x,t)&\ge -C(n,\e,a,b)+\int_{B_x(\e r)}H(x,y,t)u(y)\,dy\\
&\ge -C(n,\e,a,b)+
\frac{C_2}{V_x(\sqrt{ t})}\int_{B_x(\e \sqrt t)}\exp(-\frac{r^2(x,y)}{5t}){\Cal
B}(y)\, dy \\
&\ge -C(n,\e,a,b)+ \frac{C_3V_x(\e \sqrt t)}{V_x(\sqrt{t})} \inf_{B_x(\e r)}u \\
&\ge -C(n,\e,a,b)+C_4\inf_{B_x(\e r)}u
\endsplit
$$
for some positive constants  $C_2-C_3$  depending only on $n$ and $\e$.
The proof of the corollary is completed.
\enddemo

Suppose $u$ is Lipschitz, so that $|u(x)-u(y)|\le \beta r(x,y)$, then  $v$ is defined for all $t$.
 We have the following.
\proclaim{Lemma 1.4}    Suppose $u$ is Lipschitz so that $|u(x)-u(y)|\le \beta r(x,y)$ for all $x,\ y\in M$ and let $v$ be the solution of the heat equation with initial value $u$ obtained in Lemma 1.2.  Then for all $t>0$,
$$
\sup_{x\in M }|\nabla v(x,t)|\le  \beta.
$$
\endproclaim
\demo{Proof} By \cite{GW 1, Proposition 2.1}, for any $i>0$, there is a
smooth function $u_i$ such that $\sup_{ M}|\nabla u_i|\le
\beta+i^{-1}$ and $\sup_{M}|u_i-u|\le i^{-1}$. By Lemma 1.2, we
can solve the initial value problem for the heat equation with
initial value $u_i$. Denote the solution by $v_i$, which is
defined for all $t$. Moreover,
$$
|v-v_i|(x,t) \le \int_{M}H(x,y,t)|u(y)-u_i(y)|dy\le i^{-1}.
$$
In particular, for $x\in M$ and $t>0$, after passing to a subsequence if necessary,
$$
\lim_{i\to\infty}|\nabla v_i|(x,t)=|\nabla v|(x,t).\tag1.10
$$
However, using a more general version of  \cite{LT, Proposition 2.4}, see Lemma 1.5 below, we have
$$
\sup_{M}|\nabla v_i|(x,t)\le \sup_{M}|\nabla u_i|\le \beta+i^{-1}.\tag1.11
$$
The lemma follows from (1.10) and (1.11).
\enddemo

\proclaim{Lemma 1.5} Let $M$ be a complete noncompact Riemannian manifold with nonnegative Ricci curvature. Let $u$ be a smooth function on $M$ with bounded gradient and let $v$ be the solution of the heat equation initial value $u$. Then
for any $t>0$
$$
\sup_M|\nabla v(\cdot,t)|\le \sup_M|\nabla u|.
$$
\endproclaim
\demo{Proof} For any $T>0$, by Lemma 1.3, since $|u|$ is of linear growth, we have
$$
|v(x,t)|\le C_1\lf(r(x)+1\ri)
$$
for some $C_1$ for all $(x,t)\in M\times[0,T]$. On the other hand, Using the fact that $\heatt v^2=2|\nabla
v|^2$, and  using a suitable cut off function, one can obtain
$$
\split \int_0^{T}\int_{B_o(r)}|\nabla v|^2dxdt&\le
C_{18}\lf[r^{-2}\int_0^T\int_{B_o(2r)}v^2dxdt+\int_{B_o(2r)}u^2dx\ri]
\endsplit
$$
and so
$$
\int_0^T\int_M \exp\lf(-r^2(x)\ri)|\nabla v|^2 dxdt<\infty.
$$
Combining with the fact that $|\nabla v|$ is a subsolution of the
heat equation, the lemma follows from the maximum principle in
\cite{KL, Theorem 1.2 of NT 1}.
\enddemo

\proclaim{Lemma 1.6} Let $M^n$ be a complete  Riemannian manifold
with nonnegative Ricci curvature. Assume that $g(x)$ is a smooth
function satisfying
$$
\Delta g\ge f \tag 1.12
$$
for some continuous function $f(x)$. Assume that $f\ge -a$ for
some $a\ge 0$ and there exists a monotone nondecreasing function
$k(r)$ such that
$$
g(x)\le k(r(x)). \tag 1.13
$$
Then
$$
\int_0^{\frac12r} s\lf(\aint_{B_o(s)} f_{+}(y) \, dy\ri)\,ds \le C(n)\lf( k(5r)-g(o)+
a\, r^2\ri)  \tag 1.14
$$
where  $\fp=\max\{f,0\}$.
In particular,
$$
r^2\aint_{B_o(r)}\fp \le C(n) \lf(k(10r)-g(o)+ar^2\ri).\tag 1.15
$$
\endproclaim
\demo{Proof} Let $M_1=M\times\R$ and let $g_1(x,t)=g(x)+\frac12 at^2$ for $(x,t)\in M\times\R$. Then $\D_{M_1}g_1\ge0$. By  Theorem 2.1 of [NST 1], we have
$$
C(n)\int_0^rs\lf(\aint_{B_{o_1}(s)}\D_{M_1}g_1\ri)\,ds\le   \sup_{B_{o_1}(5r)}g-g_1(o_1)\le \sup_{B_o(5r)}g+\frac{25}2ar^2-g(0)
$$
for some positive constant $C(n)$ depending only on $n$. Here
$o_1=(o,0)$ and $B_{o_1}(s)$ is the geodesic ball in $M_1$ with
center at $o_1$ and radius $s$. The lemma follows from the fact
that $\D_{M_1}g_1=\D g+a\ge \fp$ and Lemma 1.1 in \cite{NST 1}.
\enddemo


\input amstex
\documentstyle{amsppt}
\magnification=1200 \hsize=13.8cm \catcode`\@=11
\def\NoLogo{\let\logo@\empty}
\catcode`\@=\active \NoLogo

\def\heat{\lf(\frac{\p}{\p t}-\Delta\ri)}
\def \b {\beta}
\def\i{\sqrt{-1}}
\def\Ric{\text{Ric}}
\def\lf{\left}
\def\ri{\right}
\def\bbar{\bar \beta}
\def\a{\alpha}
\def\ol{\overline}
\def\g{\gamma}
\def\e{\epsilon}
\def\p{\partial}
\def\delbar{{\bar\delta}}
\def\ddbar{\partial\bar\partial}
\def\dbar{\bar\partial}

\def\C{\Bbb C}
\def\R{\Bbb R}
\def\tx{\tilde x}
\def\vp{\varphi}

\def\tN{\tilde\nabla}

\def\dbar{\bar\partial}
\def\ba{{\bar\alpha}}
\def\bb{{\bar\beta}}
\def\abb{{\alpha\bar\beta}}

\def\i{\sqrt {-1}}

\def \D {\Delta}
\def\aint{\frac{\ \ }{\ \ }{\hskip -0.4cm}\int}
\documentstyle{amsppt}
\vsize=19.0 cm

\subheading{\S2 A maximum principle for tensors}

In this section,  we always assume that $M^m$ is a complete
noncompact K\"ahler manifold of complex dimension $m$ (real
dimension $n=2m$). We denote the  K\"ahler metric by $g_\abb$. We
want to establish a maximum principle for   Hermitian symmetric
$(1,1)$ tensor $\eta$ satisfying the complex Lichnerowicz heat
equation:
$$
\left(\frac{\p}{\p t}-\D \right)\eta_{\g\delbar}= R_{\beta
\bar{\a}\g\delbar}\eta_{\a\bbar}-
\frac{1}{2}\left(R_{\g\bar{p}}\eta_{p\delbar}+
R_{p\delbar}\eta_{\g\bar{p}}\right).\tag2.1
$$

Assume $\eta(x,t)$ is defined on $M\times[0,T]$ for some $T>0$. We also assume that there exists a constant $a>0$ such that
$$
\int_M \|\eta\|(x,0)\exp\lf({-a r^2(x) }\ri)\, dx <\infty \tag 2.2
$$
and
$$
\liminf_{r\to\infty}\int_0^T\int_{B_o(r)}
\|\eta\|^2(x,t)\exp\lf({-ar^2(x)}\ri)\, dx\, dt <\infty. \tag 2.3
$$
Here $\|\eta\|$ is the norm of $\eta_{\abb}$ with respect to the
K\"ahler metric. By (2.2), we have
$$
\int_{B_o(r)}||\eta||(x,0)\,dx\le \exp(ar^2)\cdot \Cal S\tag2.4
$$
where $\Cal S=\int_M \|\eta\|(x,0)\exp\lf({-a r^2(x) }\ri)\, dx$.

In the following, we always arrange the eigenvalues of $\eta$ at a point in the ascending order.

Before we state our result, let us first fix some notations. Let
$\varphi:[0,\infty)\to[0,1]$ be a smooth function so that
$\varphi\equiv1$ on $[0,1]$ and $\varphi\equiv0$ on
$[2,\infty)$. For any $x_0\in M$ and $R>0$, let $\varphi_{x_0,R}$
be the function defined by
$$
\varphi_{x_0,R}(x)=\varphi\lf(\frac{r(x,x_0)}{R}\ri).
$$
Let $f_{x_0,R}$ be the solution of
$$
\heat f=-f
$$
  with initial value $\varphi_{x_0,R}$. Then $f_{x_0,R}$ is defined for all $t$ and is positive and bounded for $t>0$. In fact
$$
f_{x_0,R}(x,t)=e^{-t}\cdot \int_M H(x,y,t)\vp_{x_0,R}(y)dy.
$$

We shall establish the following maximum principle.

\proclaim{Theorem 2.1} Let $M$ be a complete noncompact K\"ahler manifold
with nonnegative holomorphic bisectional curvature. Let $\eta
(x,t)$ be a Hermitian symmetric (1,1) tensor satisfying (2.1) on
$M\times [0,T]$ with $0<T<\frac{1}{40a}$  such that $||\eta||$ satisfies (2.2) and (2.3).  Suppose at $t=0$, $\eta_\abb\ge
-bg_\abb$ for some constant $b\ge0$. Then there exists $0<T_0<T$
depending only on $T$ and $a$ so that the following are true.
 \roster
\item"{(i)}" $\eta_\abb(x,t)\ge -bg_\abb(x)$ for all $(x,t)\in
M\times[0,T_0]$.
\item"{(ii)}" For any $T_0>t'\ge0$, suppose there
is a point $x'$ in $M^m$ and there exist constants $\nu>0$ and $R>0$ such that the sum of the first
$k$   eigenvalues $\lambda_1,\dots,\lambda_k$ of $\eta_\abb$
satisfies
$$
\lambda_1+\dots+\lambda_k\ge -kb+\nu k\varphi_{x',R}
$$
for all $x$ at time $t'$. Then for
all $t>t'$ and for all $x\in M$, the sum of the first $k$
eigenvalues of $\eta_\abb(x,t)$ satisfies
$$
\lambda_1+\dots+\lambda_k\ge  -kb+\nu kf_{x',R}(x,t-t').
$$
\endroster
\endproclaim
\proclaim{Remark 2.1} It is well-known that the maximum principle
for the heat equation is not true in general. The assumption of
(2.3) type is the weakest and has been appeared for the scalar
heat equation in \cite{KL}, \cite{NT 1}. From this consideration,
(2.3) is necessary. Also (2.2) in a sense ensures the solvability
of the Cauchy problem of (2.1). Therefore, it is a reasonable
assumption.

\endproclaim

To prove the theorem, we begin with some lemmas. By Lemma 1.2 and (2.4), if we let
$$
h(x,t)=\int_M H(x,y,t)\|\eta\|(y,0)\,dy,
$$
then $h(x,t)$ is a solution of the heat equation defined on
$M\times[0,\frac{1}{40a}]$ with initial value $||\eta||$. In the
following $T_0$ always denotes a constant depending only on $a$
and satisfying $T>T_0>0$. However, it may vary from line to line.

\proclaim{Lemma 2.2} Let $M^m$ be a complete noncompact K\"ahler manifold
with nonnegative holomorphic bisectional curvature. Let $\eta$ be
a Hermitian symmetric (1,1) tensor satisfying (2.1)  on $M\times
[0, T]$. Then $\|\eta\|(x,t)$ is a sub-solution of the heat
equation. Moreover, if $\eta$ also satisfies (2.2) and (2.3), then
there exists $T>T_0>0$ depending only on $a$   such that
$\|\eta\|(x,t)\le h(x,t)$ in $M\times[0,T_0]$.
\endproclaim
\demo{Proof} The first part is direct calculation. In fact using
(2.1) one has
$$
\split \lf(\D-\frac{\p}{\p t}\ri)\|\eta\|^2& = \|\eta_{\abb
s}\|^2+\|\eta_{\abb \bar{s}}\|^2+
2R_{\a\bar{p}}\eta_{p\bar{\delta}}\eta_{\delta
\bar{\a}}-2R_{\abb p\bar{q}}\eta_{\bar{p}q}\eta_{\bar{\a}\delta}\\
& \ge \|\eta_{\abb s}\|^2+\|\eta_{\abb \bar{s}}\|^2.
\endsplit
$$
Combining with the observation
$$
2|\nabla \|\eta\||^2 \le \|\eta_{\abb s}\|^2+\|\eta_{\abb
\bar{s}}\|^2
$$
we have that $\lf(\D-\frac{\p}{\p t}\ri) \|\eta\|\ge 0$.

Since $F=||\eta||-h$ is also a subsolution of the heat equation,
the second conclusion follows from (2.3) and the proof of Theorem 1.2 of [NT
1] because the positive part of $F$ is less than or equal to
$||\eta||$.
\enddemo

Let  $A_o(r_1, r_2)$ denote the annulus $B_o(r_2)\setminus
B_o(r_1)$. For any $R>0$, let $\sigma_R$ be a cut-off function which
is $1$ on $A_o(\frac{R}{4}, 4R)$ and $0$ outside $A_o(\frac{R}{8},
8R)$. We define
$$
h_R(x,t)=\int_M H(x,y,t)\sigma_R(y)||\eta||(y,0)dy.
$$
Then $h_R$ satisfies the heat equation with initial data
$\sigma_R||\eta||$.

\proclaim{Lemma 2.3} Under the assumption
 (2.2) on $\eta$,  there exists $T_0>0$ depending only on $a$ such that the following are true.
\roster
\item"{(i)}" There exists a function $\tau=\tau(r)>0$ with
$\lim_{r\to \infty}\tau(r)=0$ such that for all $R\ge \max\{ \sqrt{T_0}, 1\}$ and for all  $(x,t)\in
A_o(\frac R2, 2R)\times[0,T_0]$,
$$
h(x,t)\le h_R(x,t)+\tau(R).
$$
\item"{(ii)}" For any $r>0$,
$$
\lim_{R\to\infty}\sup_{B_o(r)\times [0,T_0]}h_R=0.
$$
\endroster
\endproclaim

\demo{Proof} Note that  $h$ is defined on $M\times[0,\frac1{40a}]$, the first condition on $T_0$ is that $T_0<\frac1{40a}$.

(i) Suppose   $R^2\ge \max\{T_0 , 1\}$, where $T_0$ will be chosen later. For $0<t<T_0$, by the definition of $h$ and $h_R$, we have
$$
h(x,t)\le h_R(x,t)+\int_{M\setminus
B_o(4R)}H(x,y,t)||\eta||(y,0)dy+\int_{B_o(\frac
R4)}H(x,y,t)||\eta||(y,0)dy.\tag2.5
$$
For $x\in A_o(\frac R2,2R)$ and $y\in B_o(\frac R4)$, $r(x,y)\ge
\frac R4$. Moreover, $V_x(R)\ge C(m)V_o(R)$ for some constant
$C(m)>0$ by the volume comparison. Hence if $x\in A_o(\frac
R2,2R)$, using (2.4), \cite{LY, p.176} and volume comparison we
have
$$
\split
\int_{B_o(\frac R4)}H(x,y,t)||\eta||(y,0)dy&\le
\sup_{y\in B_o(\frac R4)}H(x,y,t)\int_{B_o(\frac R4)}||\eta||(y,0)\, dy\\
&\le \frac{C_1}{V_x(\sqrt t)}\sup_{y\in B_o(\frac R4)}\exp\lf(-\frac{r^2(x,y)}{5t}\ri)\int_{B_o(\frac R4)}||\eta||(y,0)\, dy\\
&\le \frac{C_2}{V_o(R)}\lf(\frac{R}{\sqrt
t}\ri)^{2m}\exp\lf(-\frac{R^2}{100t}+\frac{aR^2}{16}\ri)\cdot \Cal S,
\endsplit\tag2.6
$$
for some constants $C_1$, $C_2$ depending only on $m$. On the other hand,  since $x\in B_o(2R)$, then as in the proof of (1.3), if $T_0$ is small enough depending only on $a$, we have for $0<t\le T_0$
$$
\split \int_{M\setminus B_o(4R)}H(x,y,t)||\eta||(y,0)\,dy&\le
C_3t^{-1}  \lf[\int_{4R}^\infty  \exp(-\frac{s^2}{40t})s\aint_{  B_o(s)}||\eta||(y,0)\,ds\ri]\\
&\le C_4\Cal S\int_{4R}^\infty \exp(-\frac{s^2}{40t}+as^2)\,d\lf(\frac{s^2}{t}\ri)\\
&\le C_4\Cal S\int_{\frac{16R^2}{T_0}}^\infty\exp(-\frac1{80}\tau)\,d\tau
\endsplit\tag2.7
$$
where $C_3-C_4$ are constants depending only on $m$, provided $T_0$ is small enough. Here we have used (2.4) and the fact $t\le T_0$.
 From (2.5)--(2.7), (i) follows.

 (ii) For $r>0$ fixed, if $T_0$ is small enough depending only on $a$ and
 $(x,t)\in B_o(r)\times(0,T_0]$, for $R>>r$, as before we have
$$
\split
h_R(x,t)&\le\int_{M\setminus B_o(\frac R8)}H(x,y,t)||\eta||(y,0)dy\\
&\le C_8\int_{\frac R8}^\infty \exp(-\frac{s^2}{100t})\aint_{B_o(s)}||\eta||(y,0) ds\\
&\le C_9\Cal S\int_{\frac R8}^\infty\exp(-\frac{s^2}{100t}+as^2)
ds,
\endsplit
$$
for some constants $C_8$ and $C_9$ depending only on $m$.  From this (ii) follows.
\enddemo
Now we are ready to prove Theorem 2.1.

\demo{Proof of Theorem 2.1} We only prove (b)   by assuming (a) is
true.  The proof of (a) is similar. Without loss of generality we
assume that $t'=0$. Let $T\ge T_0>0$ be small enough  so that
Lemmas 2.2, 2.3 are true. $T_0$ depends only on $T$ and $a$. By Lemma 1.2 and Corollary 1.4, we can find a solution $\phi(x,t)$,
$$
\heat \phi=\phi \tag2.8
$$
such that $\phi(x,t)\ge \exp(c(r^2(x)+1))$ for some $c>0$ for  all
$0\le t\le T$. For example, let  $\phi(x,t)=e^t\int_MH(x,y,t)h(y)dy$ with $h(y)\ge \exp(c'r^2)$ for some $c'>0$.

Assume at $t=0$, there exist $x_0\in M$,  $\nu>0$
and $R>0$ such that the first  $k$   eigenvalues
$\lambda_1,\dots,\lambda_k$ of $\eta_\abb$ satisfy
$$
\lambda_1+\dots+\lambda_k\ge -kb+\nu k\varphi_{x_0,R}
$$
for all $x$ at time $t=0$. For simplicity, let us assume that $\nu=1$.

By Lemma 2.2 and Lemma 2.3(i)
$$
||\eta_\abb||(x,t)\le  h(x,t)\le h_R(x,t)+\tau(R)\tag2.9
$$
for all $(x,t)\in \p B_o(R)\times[0,T_0]$,  where $\tau(R)>0$ is a
constant depending only on $R$ and $\tau(R)\to0$ as $R\to\infty$.

Let $\e>0$, for any $R>0$,  define $\psi=-f+\e\phi+h_R+\tau(R)+b$,
where $h_R$ is the function defined above and $f(x,t)=f_{x_0,R}(x,t)$. Let
$(\eta_R)_\abb=\eta_\abb+\psi g_\abb$, where $g_\abb$ is the
metric tensor of $M$. Then at $t=0$, at each point the sum of the
first $k$  eigenvalues of $\eta_R $ is positive. We want to prove
that for any $T_0\ge  t>0$ and $R>0$, the sum of the first
$k$ eigenvalues of $\eta_R$ in $B_o(R)\times [0,T_0]$ is positive, provided $R$ is large enough.

Let $R$ be large enough so that $\e\phi -f>0$ outside $B_o(R)$.
Then by the definition of $\psi$ and (2.9), $(\eta_R)_\abb$ is
positive definite on $\p B_o(R)\times[0,T_0]\cup
B_o(R)\times\{0\}$ and hence it is positive definite in a
neighborhood of this set. Suppose there exists $(x,t)\in \overline
B_o(R)\times[0,T_0]$ such that the sum of the first $k$
eigenvalues of $(\eta_R)_\abb$ is negative, then there exists
$0<t_1\le T_0$ and a point $x_1\in \overline B_o(R)$ such that  the
sum of the first $k$ eigenvalues of $\eta_R$ at $x_1$ at time
$t_1$ is zero but the sum of the first $k$ eigenvalues of $\eta_R$
at any point $(x,t)\in B_0(R)\times[0,t_1)$ is positive.

Let us fix the notations. Suppose $v_1,\dots,v_m$ are unit
eigenvectors of $\eta_R$ at $(x_1,t_1)$, with eigenvalues
$\lambda_1\le \lambda_2\le\dots\le \lambda_m$.   We may choose
normal coordinates at $x_1$   such that $v_j=\frac{\p}{\p z^j}$ at
$x_1$. In particular, if we write $v_j=v_j^\a\frac{\p}{\p z^\a}$,
we have $v_j^\a=\delta_{\a j}$ at $x_1$.  Note  that the sum of
the first $k$ eigenvalues of a Hermitian form is   the infimum of
the traces of the form  restricted to
  $k$-dimensional subspaces. Therefore $\sum_{\a, \beta =1}^{k}
\lf(g^{\abb}(\eta_R)_{\abb}\ri)\ge 0$ for all $(x,t)$ with $t\le
t_1$ and equals to zero at $(x_1, t_1)$. Since $\eta_R$ is
positive definite in a neighborhood of $\p B_o(R)\times[0,T_0]$, we
conclude that $x_1$ is an interior point on $B_o(R)$.

Hence at $(x_1,t_1)$, we have
$$
0\ge \heat\lf(\sum_{ \a,\beta =1 }^k  (\eta_R)_\abb
g^{\abb}\ri).\tag2.10
$$
From now on repeated indices mean summation from $1$ to $m$ if
there is no specification. Now
$$
\frac{\p}{\p t}\lf(\sum_{\a,\beta =1}^k  (\eta_R)_\abb
g^{\abb}\ri) =\sum_{\a,\beta=1}^k\lf(\frac{\p}{\p t}(\eta_R)_\abb
\ri)g^{\abb}. \tag2.11
$$
 Also at $(x_1,t_1)$, we have
$$\Delta \lf(\sum_{\a,\beta=1}^k  (\eta_R)_\abb g^{\abb}\ri)
= \sum_{\a,\beta=1}^k \lf(\Delta(\eta_R)_\abb\ri) g^{\abb}\tag2.12
$$
By   (2.10)--(2.12) and (2.1), at $(x_1,t_1)$ we have,
$$
\split
0&\ge \sum_{\a,\beta=1}^k \bigg[R_{\delta\bar\g\a\bb}\lf(\eta_{\g\delbar}+\psi g_{\g\delbar}\ri)-
\frac12R_{\a\bar p}\lf(\eta_{p\bb}+\psi g_{p\bb}\ri) -\frac12R_{p\bb }\lf(\eta_{\a\bar p }+\psi g_{\a\bar p }\ri)\bigg]g^{\abb}\\
&+\sum_{\a,\beta=1}^k\bigg(\lf[\heat\psi\ri]g_\abb
-R_{\delta\bar\g\a\bb}\psi g_{\g\delbar}+\frac12\psi R_{\a\bar
p}g_{p\bb} +\frac12\psi R_{\a\bar p }g_{p\bb}\bigg)g^{\abb}
.\endsplit\tag2.13
$$
Since at $(x_1,t_1)$, $\eta$ has eigenvectors $v_p=\frac{\p}{\p
z^p}$, for $1\le p\le m$, with eigenvalue $\lambda_p$
$$
\split \sum_{\a,\beta=1}^k
\bigg[R_{\delta\bar\g\a\bb}&\lf(\eta_{\g\delbar}+\psi g_{\g\delbar}\ri)-\frac12R_{\a\bar p}
\lf(\eta_{p\bb}+\psi g_{p\bb}\ri) -\frac12R_{p\bb }\lf(\eta_{\a\bar p }+\psi g_{\a\bar p }\ri)\bigg]g^{\abb}\\
&=\sum_{\a=1}^k\sum_{\gamma=1}^m R_{\gamma\bar\g \a\bar\a}\lambda_\g-\sum_{\a=1}^k R_{\a\bar\a}\lambda_\a\\
&=\sum_{\a=1}^k\sum_{\gamma=1}^m R_{\gamma\bar\g \a\bar\a}\lambda_\g-\sum_{\a=1}^k \sum_{\gamma=1}^mR_{\g\bar\g \a\bar\a}\lambda_\a\\
&=\sum_{\a=1}^k\sum_{\g=k+1}^m\lambda_\g R_{\g\bar \g \a\bar\a}-\sum_{j=1}^k\sum_{\g=k+1}^mR_{\g\bar \g \a\bar\a}\lambda_\a\\
&=\sum_{\a=1}^k\sum_{\g=k+1}^m R_{\g\bar \g \a\bar\a}(\lambda_\g-\lambda_\a)\\
&\ge 0
\endsplit\tag2.14
$$
where we have used that fact that $M$ has nonnegative bisectional
curvature, and $\lambda_\g\ge \lambda_\a$ for $\g\ge \a$. Also by
(2.8), the definition of $f$ and the fact that $\heat h_R=0$,  we have
$$
\lf[\heat\psi\ri]=f+\e\phi>0. \tag2.15
$$
Moreover
$$
 \lf(-R_{\delta\bar\g\a\bb}\psi g_{\g\delbar}+\frac12\psi R_{\a\bar p}g_{p\bb} +
 \frac12\psi R_{\a\bar p }g_{p\bb}\ri)g^{\abb}=0.\tag 2.16
$$
From (2.13)--(2.16), we have a contradiction. Hence the sum of the
first $k$ eigenvalues of $\eta_R$ is nonnegative for all $(x,t)\in
B_o(R)\times(0,T_0]$. In particular, if $r>0$ is fixed then the sum of
  the first $k$ eigenvalues of $\eta_R$ is nonnegative for all
$(x,t)\in B_o(r)\times(0,T_0]$. Let $R\to\infty$, using  Lemma
2.3, we conclude that the sum of the first $k$   eigenvalues of
$$
\eta_\abb(x,t)+\lf(-f(x,t)+\e\phi(x,t)+b\ri)g_\abb(x,t)
$$
is
nonnegative on $B_o(r)\times[0,T_0]$ and hence on
$M\times[0,T_0]$. Let $\e\to0$, we conclude that the  sum of the first
$k$  eigenvalues of $\eta_\abb(x,t)$ must be larger than or equal to
$-kb+kf(x,t)$, for $0<t\le T_0$ and for all $x\in M$.
\enddemo
From Theorem 2.1, applying the argument of [H] we can have the
following corollary.

\proclaim{Corollary 2.1} Let $M$ and $\eta$ be as in Theorem 2.1 with $b=0$.
That is  $\eta(x,0)\ge 0$ for all $x\in M$. Let $T_0>0$ such that the conclusions of the theorem is true. For $0<t<T_0$, let
$$ \Cal K(x,t)=\{w\in
T^{1,0}_x(M)|\  \eta_\abb(x,t) w^\a=0,\text{\rm \ for all $\b$}\}
$$
be the null space of $\eta_\abb(x,t)$. Then there exists
$0<T_1<T_0$ such that  for any $0<t<T_1$, $\Cal K(x,t)$ is a
distribution on $M$. Moreover the distribution is invariant under
parallel translations. In particular, if $M$ is simply-connected,
then $M=M_1\times M_2$ isometrically and holomorphically, where
${\Cal K}$ corresponds the tangent space of $M_1$,
$(\eta_{\abb}(x,t))>0$ on $M_2\times(0, T_1)$. Both $M_1$ and
$M_2$ are complete K\"ahler manifolds with nonnegative bisectional
curvature.
\endproclaim

\demo{Proof} By Theorem 2.1, $\eta(x,t)\ge0$ on $M\times[0,T_0)$. By Theorem 1.2(ii), we conclude that if $\dim\Cal K(x_0,t_0)\le k$ for some $x_0\in M$ and   $0\le t_0< T_0$ then $\dim\Cal K(x,t)\le k$ for all $x\in M$ and $t>t_0$. It is easy to see that there exists
$0<T_1<T$ such that $\dim \Cal K(x,t)$ is constant on
$M\times(0,T_1)$. Hence for each $0<t<T_1$, $\Cal K(x,t)$ is a
smooth distribution on $M$. It remains to prove that the
distribution is parallel for fixed $t$. We can proceed as in
\cite{H, Lemma 8.2}.

Fix $0<t_0<T_1$, let $x_0\in M$ and let $w_0\in \Cal K(x_0,t_0)$.
Let $\gamma(\tau)$ be a smooth curve from $x_0$ and let $w(\tau)$
be the vector field obtained by parallel translation along
$\gamma$. We want to prove that $w(\tau)$ is also in the null
space  $\Cal K(\gamma(\tau),t_0)$ at $\gamma(\tau)$. Since the argument is local, we may assume that one can   extend
$w$ to be a vector field in a neighborhood of $\gamma(\tau)$, and
then extend $w$ to be a vector field independent of time $t$. Now,
projecting $w$ onto $\Cal K(x,t)$, we have a vector field $v$ such
that  $v$ is in $\Cal K(x,t)$ for all $x$ in a neighborhood of
$\gamma$ and for all  $t$. The following computations are
performed in a neighborhood of $\gamma$.

Since
$$
\eta_{\abb}v^\a=0\tag2.17
$$
for all $\b$, we have
$$
\split
0&=\frac{\p}{\p t}\lf(\eta_{\abb}v^\a\ol{v^\b}\ri)\\
&=\lf(\frac{\p}{\p t} \eta_{\abb}\ri)v^\a\ol{v^\b}+\eta_{\abb}\frac{\p v^\a}{\p t}\ol{v^\b}
+\eta_{\abb} v^\a \frac{\p\ol{v^\b}}{\p t}\\
&=\lf(\frac{\p}{\p t} \eta_{\abb}\ri)v^\a\ol{v^\b}
\endsplit\tag2.18
$$
where we have used (2.17). Choosing a unitary frame  $e_s$ at   a
point $\gamma(\tau)$, we have
$$
\split
0&=\Delta\lf(\eta_{\abb}v^\a\ol{v^\b}\ri)\\
&=\frac12\lf(\nabla_s\nabla_{\bar s}+\nabla_{\bar s}\nabla_s\ri)\lf(\eta_{\abb}v^\a\ol{v^\b}\ri)\\
&=\lf(\Delta \eta_{\abb}\ri)v^\a\ol{v^\b}-\eta_{\abb}\nabla_{\bar
s}v^\a \nabla_s\ol{v^\b}- \eta_{\abb}\nabla_{  s}v^\a \nabla_{\bar
s}\ol{v^\b}
\endsplit\tag2.19
$$
where we have used (2.17) so that
$$
\lf(\nabla_{  s}\eta_{\abb}\ri)v^\a=- \eta_{\abb}\nabla_{  s}v^\a,
\ \lf(\nabla_{\bar  s}\eta_{\abb}\ri) v^\a=-
\eta_{\abb}\nabla_{\bar s}v^\a
$$
and their complex conjugates.

Combining with (2.1), (2.18), (2.19), we have
$$
0=R_{t\bar s\abb}\eta_{s\bar
t}v^\a\ol{v^\b}+2\eta_{\abb}\nabla_{\bar s}v^\a
\nabla_s\ol{v^\b}+2 \eta_{\abb}\nabla_{  s}v^\a \nabla_{\bar
s}\ol{v^\b}.\tag 2.20
$$
We may choose $e_s$ so that at a point $\eta_{s\bar
t}=a_s\delta_{st}$. Then
$$
R_{t\bar s\abb}\eta_{s\bar t}v^\a\ol{v^\b}=R_{s\bar s \abb }a_sv^\a
\ol{v^\b}=a_sR_{s\bar s v\bar v}\ge0
$$
because $a_s\ge0$  and $M$ has nonnegative bisectional curvature.
Hence (2.20) and the fact that $\eta\ge0$ imply that $\nabla_s v$ and $\nabla_{\bar s}v$ are
in the null space  $\Cal K(\gamma(\tau),t_0)$.

Since $w(\tau)$ is parallel along $\gamma(\tau)$, and
$w=v+w^\perp$, where $w^\perp$ is perpendicular to $\Cal
K(\g(\tau),t_0)$, we have
$$
0=\frac{D}{d\tau}w=\frac{D}{d\tau}v+\frac{D}{d\tau}w^\perp.
$$
Hence
$$
\frac{D}{d\tau}w^\perp=-\frac{D}{d\tau}v
$$
which is in $\Cal K$.

Now
$$
\frac{d}{d\tau}\langle w^\perp,w^\perp\rangle=\langle
\frac{D}{d\tau}w^\perp,w^\perp\rangle+\langle
w^\perp,\frac{D}{d\tau}w^\perp\rangle=0
$$
because $\frac{D}{d\tau}w^\perp$ is in $\Cal K$ and $w^\perp$ is
perpendicular to $\Cal K$. At $\gamma(0)=x_0$, $w=v_0$ and so
$w^\perp=0$ at $\gamma(0)$. Hence $w^\perp=0$ for all $\tau$ and
so $w$ is in $\Cal K$. The last statement follows from the De Rham
decomposition.
\enddemo

\proclaim{Remark 2.2} Given the  work of \cite{H}, the main
difficulty for noncompact manifolds in the proof of Corollary 2.1
is to obtain the maximum principle Theorem 2.1.  In particular, it
can be proved more easily if we assume that $\eta$ is bounded. In
\cite{C},  results similar to the corollary are obtained
independently for the case that $\eta=\Ric$ for the K\"ahler-Ricci
flow in a
 complete noncompact K\"ahler manifolds with {\it bounded} nonnegative
holomorphic bisectional curvature. However, it seems that a
maximum principle is still needed in this case.
\endproclaim


\input amstex
\documentstyle{amsppt}
\magnification=1200 \hsize=13.8cm \catcode`\@=11
\def\NoLogo{\let\logo@\empty}
\catcode`\@=\active \NoLogo

\def\heatt{\lf (\Delta-\frac{\p}{\p t}\ri)}
\def\heat{\lf(\frac{\p}{\p t}-\Delta\ri)}
\def\fp{f_+}
\def\fm{f_-}
\def \b {\beta}
\def\i{\sqrt{-1}}
\def\Ric{\text{Ric}}
\def\lf{\left}
\def\ri{\right}
\def\bbar{\bar \beta}
\def\a{\alpha}
\def\ol{\overline}
\def\g{\gamma}
\def\e{\epsilon}
\def\p{\partial}
\def\delbar{{\bar\delta}}
\def\ddbar{\partial\bar\partial}
\def\dbar{\bar\partial}

\def\C{\Bbb C}
\def\R{\Bbb R}
\def\tx{\tilde x}
\def\vp{\varphi}

\def\tN{\tilde\nabla}

\def\dbar{\bar\partial}
\def\ba{{\bar\alpha}}
\def\bb{{\bar\beta}}
\def\abb{{\alpha\bar\beta}}

\def\i{\sqrt {-1}}

\def \D {\Delta}
\def\aint{\frac{\ \ }{\ \ }{\hskip -0.4cm}\int}
\documentstyle{amsppt}
\vsize=19.0 cm

\subheading{\S3 $C^{\infty}$-approximation to continuous
plurisubharmonic functions}\vskip .2cm

In [GW 1], it was proved that on a complete noncompact K\"ahler
manifold a continuous {\it strictly} plurisubharmonic function can
be approximated uniformly by   $C^\infty$ strictly
plurisubharmonic functions. If the function is only
plurisubharmonic, then it can be approximated uniformly by
$C^\infty$ functions whose complex Hessian are close to being
nonnegative, see Lemma 3.2 below. In general,  it seems unlikely
that a continuous plurisubharmonic function can be approximated by
$C^\infty$-plurisubharmonic functions. However, in this section,
we shall show that this can be done by the solution to the heat
equation if the K\"ahler manifold  has  nonnegative holomorphic
bisectional curvature, provided  the continuous plurisubharmonic
function satisfies a mild growth condition.  Actually, we shall
prove more in Theorem 3.1.

Let $M$ be a complete noncompact K\"ahler manifold
with nonnegative holomorphic bisectional curvature. Let $u$ be a continuous
plurisubharmonic function defined on $M$ with growth rate
satisfying
$$
|u|(x)\le C \exp(ar^2(x))  \tag 3.1
$$
 for some positive constants $a$ and $C$.
 Let $v(x,t)$ be the solution to the heat equation
  on $M\times[0,\frac{1}{40a}]$ with initial value $u$, obtained by Lemma 1.2.

\proclaim{Theorem 3.1} Let $M^m$, $u$ and $v$ be as above. There exists $T_0>0$ depending only on $a$ and there exists $T_0>T_1>0$ such that the following are true.
\roster
\item"{(i)}" For $0<t\le T_0$, $v(\cdot,t)$ is a smooth plurisubharmonic function.
\item"{(ii)}" Let
$$
\Cal K(x,t)=\{w\in T^{1,0}_x(M)|\  v_\abb(x,t) w^\a=0,\text{\rm \ for all $\b$}\}
$$
be the null space of $v_\abb(x,t)$. Then for any $0<t<T_1$, $\Cal
K(x,t)$ is a  distribution on $M$. Moreover the distribution is
invariant under parallel translations. \item"{(iii)}" If the
holomorphic bisectional curvature is positive at some point, then
$v(x,t)$ is strictly plurisubharmonic for all $0<t<T_1$.
\endroster
\endproclaim
As a corollary, we have:

\proclaim{Corollary 3.1} Let $M$ be a complete noncompact K\"ahler manifold with nonnegative holomorphic bisectional curvature and let $u$ be a continuous plurisubharmonic function on $M$ satisfying (3.1). Then there exist $C^\infty$ plurisubharmonic functions $u_i$ such that $u_i$ converges to $u$ uniformly on compact subsets of $M$. If in addition, the holomorphic bisectional curvature is positive at some point, then $u_i$ can be chosen to be strictly plurisubharmonic.
\endproclaim

We shall prove Theorem 3.1 by using the results in \S2 together
with a result in \cite{GW 1}. In order to use the results in \S2,
we need to the following estimates.

\proclaim{Lemma 3.1} Let $M^m$ be a complete noncompact K\"ahler
manifold with nonnegative holomorphic bisectional curvature. Let
$u$ be a {\it smooth} function satisfying (3.1) and let $v$ be the
solution of the heat equation on $M\times[0,\frac{1}{40a}]$ with
initial value $u$, obtained in Lemma 1.2. Moreover, assume that
there exists $1\ge b\ge0$ such that
$$
u_\abb(x)\ge -bg_\abb(x)\tag3.2
$$
for all $x\in M$. Let $||\rho||(x,t)$ be the norm of $v_\abb(x,t)$. Then there exists $\frac1{40a}>T_0>0$ depending only on $a$ with the following properties.
\roster
\item"{(i)}" There exist constants $C_1$ and $C_2$, where $C_2$ depends only on $a$ such that
$$
|v(x,t)|\le C_1\exp\lf(C_2r^2(x)\ri).
$$
for all $(x,t)\in M\times[0,T_0]$.
\item"{(ii)}" There exist constants $C_3$ and $C_4$, where $C_4$ depends only on $a$ such that
$$
\aint_{B_o(r)}||\rho||(\cdot,0)\le C_3\exp\lf(C_4r^2\ri),
$$
for all $r$.
\item"{(iii)}" There exist constants $C_5$ and $C_6$, where $C_6$ depends only on $a$ such that
$$
\int_0^{T_0}\aint_{B_o(r)}||\rho||^2(x,t)dxdt\le C_5(1+T_0)\exp\lf(C_6r^2\ri),
$$
for all $r$.
\endroster
\endproclaim
\demo{Proof} In the following $T_0$ ($\frac{1}{40a}>T_0>0$) always
denote a positive constant depending only on $a$, but its exact
value may vary from line to line.

(i) By Lemma 1.3, we conclude that there exists $T_0>0$, such that if $r=r(x)\ge \sqrt{T_0}$ then
$$
|v(x,t)|\le \int_{B_x(\frac r2)}H(x,y,t)|u|(y)dy+C_7
$$
for all $(x,t)\in M\times(0,T_0]$ for some constant $C_7$ independent of $x$ and $t$. Since $u$ satisfies (3.1) and $\int_MH(x,y,t)dy=1$, it is easy to see that (i) is true.

(ii) Let $f=\Delta u=g^\abb u_\abb$ and let $f=\fp-\fm$, where
$\fp$ ($\fm$) is the positive part (negative part)  of $f$. Let
$k^+(o,s)=\aint_{B_o(s)}\fp$. By the assumption on $u_\abb$,
$f^-\le mb\le m$. Applying (1.15) of Lemma 1.6 we have that
$$
r^2k^+(o,r) \le C(n)\lf(\exp(100ar^2)-u(o)+mr^2\ri). \tag 3.3
$$
Hence we have
$$
r^2k^+(o,r)\le C_{11}\exp(100ar^2)\tag3.4
$$
for some constant $C_{11}$ independent of $r$. On the other hand,
at a point $x$, choose an normal coordinates such
that $u_\abb=\lambda_\a g_\abb$. Since $u_\abb\ge-bg_\abb$ and $b\le1$,   for any $\alpha$
$$
\split
-1&\le-b\\
&\le \lambda_\a\\
&=\D u-\sum_{\beta\neq\a}\lambda_\beta\\
&\le \fp+(m-1)b\\
&\le \fp+(m-1).
\endsplit
$$
Therefore,
$$
||\rho||(x)\le m\lf(\fp(x)+(m-1)\ri).\tag3.5
$$
 (ii) follows from (3.4) and (3.5).

(iii) By (i), there exists $\frac{1}{40a}>T_0>0$ such that for all $(x,t)\in M\times(0,T_0)$, we have
$$
|v(x,t)|\le C_{12}\exp\lf(C_{13}r^2(x)\ri)\tag3.6
$$
for some constants $C_{12}$ independent of $x$ and $t$, and
$C_{13}$ depending only on $a$. Using $\D u=f$, integrating by
parts after multiplying a suitable cut-off function, one can prove
that
$$
\split
\int_{B_o(r)}|\nabla u|^2&\le C_{14}\lf[r^{-2}\int_{B_o(2r)}u^2+\int_{B_o(2r)}|u|\,|f|\ri]\\
&\le C_{15}V_o(r)\lf[\exp(8ar^2)+\exp(4ar^2)\aint_{B_o(2r)}|f|\ri]\\
&\le C_{16}V_o(r)\exp(C_{17}ar^2)
\endsplit\tag3.7
$$
for some constants $C_{14}-C_{16}$ independent of $r$, and
$C_{17}$ depending only on $a$. Here we have used (3.1), (ii) and
the fact that $|f|\le m||\rho||$. Using the fact that $\heatt
v^2=2|\nabla v|^2$, and  multiplying a suitable cut off function,
one can obtain
$$
\split
\int_0^{T_0}\aint_{B_o(r)}|\nabla v|^2&\le C_{18}\lf[r^{-2}\int_0^T\aint_{B_o(2r)}v^2+\aint_{B_o(2r)}u^2\ri]\\
&\le C_{19}(T_0+1)\exp(C_{20}r^2)
\endsplit\tag3.8
$$
for some constants $C_{18}-C_{19}$ independent of $r$, and
$C_{20}$ depending only on $a$. Here we have used (3.1) and (3.6).
By the Bochner formula,
$$
\heatt |\nabla v|^2\ge 2|\nabla^2 v|^2.
$$
 Multiplying this inequality by a suitable cutoff function and integrating by parts, using (3.7) and (3.8) we have
$$
\split
\int_0^{T_0}\aint_{B_o(r)}|\nabla^2 v|^2\le& C_{21}\lf[\frac1{r^2}\int_0^{T_0}\aint_{B_o(2r)}|\nabla v|^2+\aint_{B_o(2r)}|\nabla u|^2\ri]\\
&\le C_{22}(T_0+1)\exp(C_{23}r^2)
\endsplit
$$
for some constants $C_{21}-C_{22}$ independent of $r$, and
constant $C_{23}$ depending only on $a$. From this, (iii) follows.
\enddemo

We are ready to prove Theorem 3.1. We need the following
approximation result of Greene-Wu \cite{GW 1, Corollary 2 to
Theorem 4.1}.

\proclaim{Lemma 3.2} {\cite{Greene-Wu}} Let $u$ be a continuous plurisubharmonic function on $M$. For
any  and $b>0$, there is a $C^\infty$ function $w$ such that
\roster \item"{(i)}" $|w-u|\le b$ on $M$; and \item"{(ii)}"
$w_\abb\ge -bg_\abb$ on $M$.
\endroster
\endproclaim

\demo{Proof of Theorem 3.1} (i) Let $u$ and $v$ be as in the
theorem. Choose  $1>\e_i>0$ such that $\e_i\to0$ as $i\to\infty$.
By Lemma 3.2, we can find $u_i$ such that
$$|u_i-u|\le \e_i\tag3.9
$$ on $M$, and
$$(u_i)_\abb\ge -\e_ig_\abb\tag3.10
$$ on $M$. Since $u$ satisfies (3.1), each $u_i$ also satisfies (3.1). Namely,
$$
|u_i|(x)\le c_i\exp(ar^2(x)) \tag3.11
$$
for some constants $c_i$ independent of $x$. By Lemma 1.2, we can
solve the heat equation with initial data $u_i$ on
$M\times[0,\frac 1{40a}]$. The solution is denoted by $v_i$.  By
Lemma 3.1, (3.10) and (3.11), there exist a constant
$\frac1{40a}>T_0>0$ depending only on $a$ such that
$$
|v|(x,t)+|v_i|(x,t)\le d_i\exp(C_1r^2(x))\tag3.12
$$
$$
\aint_{B_o(r)}||\rho_i||(\cdot,0)\le d_i\exp\lf(C_1r^2\ri)\tag3.13
$$
and
$$
\int_0^{T_0}\aint_{B_o(r)}||\rho_i||^2(x,t)dxdt\le d_i(1+T_0)\exp\lf(C_1r^2 \ri)\tag3.14
$$
for some constants $d_i$ independent of $r$ and for some constant
$C_1$ depending only on $a$, where $||\rho_i||$ is the norm of
$(v_i)_\abb$. Here and below, $\frac{1}{40a}>T_0>0$ always denotes
a constant depending only on $a$, but it may vary from place to
place.

 Since the complex Hessian $(v_i)_\abb$ satisfies the Lichnerowicz heat equation (2.1) see \cite{Lemma 2.1, NT 2}. By (3.13), (3.14) and the maximum principle Theorem 2.1(i), there exists $\frac{1}{40a}>T_0>0$ such that
$$
(v_i)_\abb(x,t)\ge-\e_ig_\abb(x),\tag3.15
$$
for all $(x,t)\in M\times[0,T_0]$.

By (3.12), we can apply the
maximum principle of \cite{KL, NT 1}, to conclude that
$$
\sup_{M\times[0,T_0)}|v-\tilde v_i|\le \e_i.
$$
 Hence passing to a subsequence if necessary $v_i$ together with their derivatives
subconverge to $v$ uniformly on compact sets on $M\times(0,T_0)$. By (3.15),
we conclude that $v_\abb(x,t)\ge0$ on $M\times(0,T_0)$.

(ii) Let $T_0$ be as in (3.15),  which is obtained in Theorem 2.1.
Let $T_0>t_0>0$. Suppose there exists a point $x_0\in M$ such that
the sum of the first $k$ eigenvalues of $v_\abb(x_0,t_0)$
satisfies
$$
\lambda_{1}+\cdots+\lambda_{k}>0,
$$
then there exists $R>0$ and $\nu>0$ independent of $i$ such that the sum of the first $k$ eigenvalues of $(v_i)_\abb(x,t_0)$ satisfies:
$$
\lambda_{i,1}+\cdots+\lambda_{i,k}>k\nu,
$$
on $B_{x_0}(2R)$. Since $(v_i)_\abb$ satisfies (3.15), the sum of the first $k$ eigenvalues of $(v_i)_\abb$ satisfies:
$$
\lambda_{i,1}+\cdots+\lambda_{i,k}>-k\e_i+k\nu\vp_{x_0,R}
$$
at every point $x\in M$ at time $t_0$, where $\vp_{x_0,R}$ is the nonnegative function as in Theorem 2.1. By Theorem 2.1(ii), for $T_0>t>t_0$,
the sum of the first $k$ eigenvalues of $(v_i)_\abb$ at $(x,t)$ satisfies:
$$
\lambda_{i,1}+\cdots+\lambda_{i,k}\ge-k\e_i+k\nu f_{x_0,R}(x,t-t_0).
$$
where $f_{x_0,R}$ is the function defined in Theorem 2.1. Note that $f_{x_0,R}(x,s)>0$ if $s>0$. Let $i\to\infty$, we conclude that the sum of the first $k$ eigenvalues of $v_\abb$ at $(x,t)$ satisfies
$$
\lambda_{1}+\cdots+\lambda_{k}\ge k\nu f_{x_0,R}(x,t-t_0).
$$
Hence we have proved that if there exists a point $x_0\in M$ such
that the sum of the first $k$ eigenvalues of $v_\abb(x_0,t_0)$ is
positive, then for all $x\in M$ and $t>t_0$, the sum of the first
$k$ eigenvalues of $v_\abb(x,t)$ is also positive. One can then
proceed as in the proof of Corollary 2.1 to conclude that (ii) is true.

(iii) Suppose $M$ has positive holomorphic bisectional curvature
at $x_0$. For $0<t<T_0$, suppose $\dim \Cal K(x_0,t)>0$. By (ii),
locally $M$ can be  splitted  isometrically as a nontrivial
product of two K\"ahler manifold with nonnegative holomorphic
bisectional curvature. This is impossible. Hence $\dim\Cal
K(x_0,t)=0$ and so $\dim\Cal K(x,t)=0$ for all $x$. This implies
that $v(\cdot,t)$ is strictly plurisubharmonic. The proof of the
theorem is completed.
\enddemo

Using Theorem 3.1, we shall prove the following Liouville theorem
which will be used to prove a splitting theorem as well as a gap
theorem in section 4 and 6.

\proclaim{Theorem 3.2} Let $M$ be a complete K\"ahler manifold
with nonnegative holomorphic bisectional curvature. Let $u$ be a continuous
plurisubharmonic function on $M$. Suppose that
$$
\limsup_{x\to \infty} \frac{u(x)}{\log r(x)}=0. \tag 3.16
$$
Then $u$ must be a constant.
\endproclaim

To prove the theorem we need the following lemma.

\proclaim{Lemma 3.3} (\cite{N, Proposition 4.1}) Let $M^m$ be a
complete noncompact K\"ahler manifold of complex dimension $m$,
with nonnegative Ricci curvature. Let $u(x)$ be a plurisubharmonic
function on $M$ satisfying (3.16). Then $(\p\dbar u)^m =0$
\endproclaim

Proposition 4.1 stated in [N] is under the assumption that $M$ is
nonparabolic. However, the proof without any changes also works
for general complete K\"ahler manifolds with nonnegative Ricci
curvature.

\demo{Proof of Theorem 3.2} Let $M$ and $u$ satisfy the conditions
in Theorem 3.2. Let $\tilde M$ be the universal cover of $M$, then
the distance function in $\tilde M$ dominates the distance
function in $M$. Hence $\tilde M$ and the lift $\tilde u$ of $u$
also satisfy the conditions in the theorem. Therefore, we may
assume that $M$ is simply connected.

First we let $u_c=\max\{ u, c\}$. By the assumption (3.16) it is
easy to see that $u_c$ satisfying (3.1) and $u_c$ is
plurisubharmonic. Therefore, we can solve the heat equation with
$u_c(x)$ as the initial data. Denote the solution by $v_c$ on
$M\times [0, T_0]$. By adding a constant we can also assume that
$u_c(x)\ge 0$. Applying Theorem 3.1(i) to $v_c(x,t)$ we conclude that
$v_c(x,t)$ is plurisubharmonic. By Theorem 3.1(ii),  for any $t_0>0$ small enough,   $M=M_1\times M_2$ isometrically and
holomorphically such that $(v_c)_\abb$ is zero when restricted on
$M_1$ and $(v_c)_\abb$ is positive everywhere when restricted on
$M_2$ by  the De Rham decomposition (Cf. Theorem 8.1, page
172 of [K-N]).    By Corollary 1.4, we still have
$$
\limsup_{x\to\infty}\frac{v_c(x,t_0)}{\log r(x)}=0.\tag3.17
$$
Hence when restricted on $M_2$, (3.17) is still true. This
contradicts Lemma 3.3 and the fact that $(v_c)_\abb$ is positive
when restricted on $M_2$, unless $M=M_1$. Hence
$(v_c)_\abb(x,t_0)\equiv0$ on $M$ for all $0<t_0$ small enough.  By the gradient estimate of Cheng-Yau \cite{C-Y} and (3.17)  we can conclude that  that $v_c(x,t_0)$ is a constant, provided $t_0>0$ is small enough. Hence
$u_c$ is  a constant. Since $c$ is arbitrary, it shows that $u(x)$
is also a constant.
\enddemo
\input amstex
\documentstyle{amsppt}
\magnification=1200 \hsize=13.8cm \catcode`\@=11
\def\NoLogo{\let\logo@\empty}
\catcode`\@=\active \NoLogo
\def\tM{\tilde M}
\def\vgk{\text{\bf VG}_k}
\def\heatt{\lf (\Delta-\frac{\p}{\p t}\ri)}
\def\heat{\lf(\frac{\p}{\p t}-\Delta\ri)}
\def\fp{f_+}
\def\fm{f_-}
\def \b {\beta}
\def\i{\sqrt{-1}}
\def\Ric{\text{Ric}}
\def\lf{\left}
\def\ri{\right}
\def\bbar{\bar \beta}
\def\a{\alpha}
\def\ol{\overline}
\def\g{\gamma}
\def\e{\epsilon}
\def\p{\partial}
\def\delbar{{\bar\delta}}
\def\ddbar{\partial\bar\partial}
\def\dbar{\bar\partial}

\def\C{\Bbb C}
\def\R{\Bbb R}
\def\tx{\tilde x}
\def\vp{\varphi}

\def\tN{\tilde\nabla}

\def\dbar{\bar\partial}
\def\ba{{\bar\alpha}}
\def\bb{{\bar\beta}}
\def\bg{{\bar\gamma}}
\def\abb{{\alpha\bar\beta}}

\def\i{\sqrt {-1}}

\def \D {\Delta}
\def\aint{\frac{\ \ }{\ \ }{\hskip -0.4cm}\int}
\documentstyle{amsppt}
\vsize=19.0 cm

\subheading{\S4 Structure of  non-negatively curved K\"ahler
manifolds I}

In this and the next section, we shall apply the results in the previous two sections to study the structure of complete noncompact K\"ahler manifolds with nonnegative sectional or  holomorphic bisectional curvature. Let us  begin with some lemmas.

Let $M^m$ be a complete noncompact K\"ahler manifold with
nonnegative holomorphic bisectional curvature. Recall the
definition of the Busemann function at a point $o\in M$, see
\cite{CG 2}. Let $\gamma$ be a ray from $o$ parametrized by arc length.
Then the  Busemann ${\Cal B}_\gamma(x)$ is defined as
$$
\Cal B_\gamma(x)=\lim_{s\to\infty}\lf(s-d(x,\gamma(s)\ri)
$$
and the Busemann function $\Cal B$ is defined as
$$
\Cal B(x)=\sup_\gamma \Cal B_\gamma(x)
$$
where the supremum is taken over all rays from $o$. It is well-known that $\Cal B$ is Lipschitz with Lipschitz constant 1. Since $M$ has nonnegative holomorphic bisectional curvature, by the result of  Wu \cite{W, p.58}, $\Cal B$ is a continuous plurisubharmonic function on $M$.  Let $v(x,t)$ be the solution of the heat equation with initial value $\Cal B$. Then $v$ is defined for all $t$. We collect some facts in the follow lemma for easy reference.

\proclaim{Lemma 4.1} With the above assumptions and notations, the following are true.
\roster
\item"{(i)}" For any $t>0$, $v(\cdot,t)$ is a smooth  plurisubharmonic function. If the holomorphic bisectional curvature of $M$ is positive at some point, then $v(\cdot,t)$ is strictly plurisubharmonic.
\item"{(ii)}" For any $t>0$,
$$
\sup_M|\nabla v(\cdot,t)|\le 1.
$$
\item"{(iii)}"  For any $t>0$, $v(\cdot,t)$ grows linearly  when restricted on a ray from $o$. If in addition,   $\Cal B$ is an exhaustion function of $M$, then $v(\cdot,t)$ is also an exhaustion function of $M$ for all $t>0$.
\item"{(iv)}" There exists $T_0>0$, such that for any $0<t<T_0$, the null space $\Cal K(x,t)\subset T^{(1,0)}_x(M)$  of $v_\abb(x,t)$ is a parallel distribution on $M$.
\endroster
\endproclaim
\demo{Proof} (i) and (iv) are just special cases of Theorem 3.1. Note that from the proof of Theorem 3.1, (i) is true for any $t>0$. (ii) follows from Lemma 1.4. It remains to prove (iii). Let $\gamma$ be a ray from $o$ and let $x=\gamma(r)$ where $r=r(x,o)$. Then $\Cal B(x)\ge \Cal B_\gamma(x)\ge r$. Since ${\Cal B}$ has
Lipschitz constant $1$ we know that ${\Cal B}(y)\ge \frac{1}{2}
r$, for all $y\in B_x(\frac{r}{2})$. By Corollary 1.4,  we
know that
$$
v(x,t)\ge C_1r-C_2  \tag4.1
$$
for some positive constants  $C_1$ and $C_2$   independent of $x$.  From this we conclude that $v(\cdot, t)$ grows linearly on $\gamma$. The second part of (iii) can be proved similarly.
\enddemo

Recall that  $M$ is said to satisfy ($\vgk$) for $k>0$, if there exists a constant $C>0$ such that
$$
V_o(r)\ge Cr^{k}\tag"($\vgk$)"
$$
for all $r\ge 1$.   $M$ is
said to satisfy the curvature decay condition ({\bf CD}) if there
exists a constant $C>0$  such that
$$
\aint_{B_o(r)}\Cal R\le \frac{C}{r}\tag"({\bf CD})"
$$
for all $r>0$. Finally, $M$ is said to satisfy the fast curvature
decay condition ({\bf FCD}) if there is a constant $C>0$, so that
$$
\int_0^rs\lf(\aint_{B_o(s)}\Cal R(x)dx\ri)ds\le C\log(r+2)\tag"({\bf FCD})"
$$
for all $r>0$.

  Using the ideas in \cite{CZ 2}, we can prove the following.

\proclaim{Lemma 4.2} Let $M^m$ be a complete noncompact K\"ahler
manifold.   \roster \item"{(i)}" Suppose $M$ supports a smooth
plurisubharmonic function $u$ which is strictly pluri-subharmonic
at $o$ and suppose $u$ has   bounded gradient. Then $M$ satisfies
($\text{\bf VG}_m$), where $m$ is the complex dimension of $M$. If
in addition, $M$ has nonnegative Ricci curvature and also supports a nonconstant holomorphic function
of polynomial growth then $M$ satisfies ($\text{{\bf VG}}_{a}$)
for any $a<m+1$.
 \item"{(ii)}" Suppose $M$
has nonnegative Ricci curvature and suppose $M$ supports a
strictly plurisubharmonic function $u$. If $u(x)\le C(r(x)+1)$ for
some constant $C$,  then $M$ satisfies ({\bf CD}). If  $u(x)\le
C\log (r(x)+1)$ for some constant $C$, then $M$ satisfies ({\bf
FCD}).
\endroster
\endproclaim
\demo{Proof} (i) Let $\omega$ be the K\"ahler form of $M$ which is closed. Since $\i\ddbar u\ge0$ and $\i\ddbar u>0$ at $o$, for any $r>1$, there exists a smooth cutoff function $0\le\vp\le 1$ such that $\vp\equiv1$ on  $B_o(r)$ and $\vp\equiv0$ outside $B_o(2r)$ and such that $|\nabla \vp|\le C_1/r$ for some constant $C_1$ independent of $r$ and
$$
\split
\int_{B_o(1)}\lf(\i\ddbar u\ri)^m &\le \int_{B_o(2r)}\vp^m\lf(\i\ddbar u\ri)^m \\
&=-m\int_{B_o(2r)} \vp^{m-1}\i\p\vp\wedge\dbar u\wedge \lf(\i\ddbar u\ri)^{m-1}\\
&\le \frac{mC_2}{r}\int_{B_o(2r)}\vp^{m-1}\lf(\i\ddbar u\ri)^{m-1}\wedge\omega
\endsplit
$$
for some constant $C_2$ independent of $r$, where we have used the fact that $|\nabla \vp|\le C_1/r$ and $|\nabla u|$ is bounded.
 Continuing  in this way and integrating by parts $(m-1)$ times more, we have
$$
\int_{B_o(1)}\lf(\i\ddbar u\ri)^m\le   {m!}\cdot\lf(\frac{C_2}{r}\ri)^mV_o(2r).
$$
  Since $\ddbar u>0$ at $o$, it is easy to see that $M$ satisfies ($\text{\bf VG}_m$).

If in addition, $M$ has nonnegative Ricci curvature and supports a nonconstant polynomial growth holomorphic function $f$. Let $v(x)=\log (|f|^2+1)$. Then    $v(x)\le C\log (r(x)+2)$,   and $v$ is plurisubharmonic. Moreover, $\ddbar v$ is not zero at every point outside a subvariety. Observe that
$$
r^2\aint_{B_o(r)}\D v(y)\, dy\le C_3\log (r+2)\tag4.1
$$
by Lemma 1.6 for some constant $C_3$ independent of $r$. On the other hand, using the cut-off function $\vp$
as above, we have that
$$
\split 0<\int_{B_o(1)}\lf(\sqrt{-1}\ddbar v\ri)\wedge
\lf(\sqrt{-1} \ddbar u\ri)^{m-1} &\le
\int_{B_o(2r)}\vp^m \lf(\sqrt{-1}\ddbar v\ri)\wedge\lf(\sqrt{-1}\ddbar u\ri)^{m-1}\\
&\le \frac{C_4}{r^{m-1}}\int_{B_o(2r)}\sqrt{-1}\ddbar
v\wedge\omega^{m-1}\\
&\le \frac{C_5}{r^{m-1}}\int_{B_o(2r)}\Delta v(y)\, dy
\endsplit\tag4.2
$$
for some constants $C_4-C_5$ independent of $r$. Combining (4.1) and (4.2), we have that for some positive constant
$C_6$ independent of $r$ such that,
$$
V_o(r)\ge C_6\frac{\log (r+2)}{r^{m+1}}.
$$
This concludes the proof of  (i).

(ii) Let us prove the second statement.  Our proof  is basically a
simplified   version of  \cite{CZ 2}.  Using $u$ as a weight
function, by the $L^2$ estimate and  Theorem 3.2 of [N 1],   there
exists a nontrivial holomorphic section $s$ of the canonical line
bundle $K_S$ (a $(n,0)$ form in terms of Theorem 3.2 of [N 1]) such
that $s(o)\neq0$ and
$$
\int_M \|s\|^2\exp(-C_7 u(x)) \, dx= {\Cal A}<\infty. \tag 4.3
$$
for some constant $C_7>0$. Since $u(x)\le C\log(r(x)+2)$, for some constant $C$ independent of $x$, (4.3)   implies  that
$$
\int_{B_{o}(R)}\|s\|^2(x)\, dy \le  (R+1)^{C_8}$$
for some constant $C_8$ independent of $R$.  It is well-known that  $\|s\|^2$ is subharmonic, see Lemma 4.2 of [NST 2] for example.  By the mean value inequality
of Li-Schoen \cite{LS, p.287},  we have that
$$
\|s\|^2(x)\le C(m)\aint_{B_{o}(2r(x))}\|s\|^2(y)\, dy
$$
for some constant $C(m)$ depending only on $m$. Therefore we have that
$$
\|s\|^2(x)\le  (r(x)+1)^{C_9}
$$
for some constant $C_9$ independent of $x$, and so
$$
\log\lf(||s||^2(x)+1\ri)\le C_{10}\log(r(x)+2)\tag4.4
$$
for some constant $C_{10}$ independent of $x$.   By Lemma 4.2 of [NST 2] again, for any $1>\e>0$,  we have that
$$
\Delta \log (\|s\|^2(x) +\e)\ge {\Cal
R}(x)\cdot\frac{\|s\|^2}{\|s\|^2+\e}  \tag 4.5
$$
where $\Cal R$ is the scalar curvature of $M$.
Applying  Lemma 1.6, noticing that ${\Cal
R}\cdot\frac{\|s\|^2}{\|s\|^2+\e}\ge 0$, we have that
$$
\int_0^r \sigma\lf(\aint_{B_{o}(\sigma)}{\Cal R}(x)\cdot
\frac{\|s\|^2(x)}{\|s\|^2(x)+\e}\, dx \ri)\, d\sigma \le C_{11} \log(r+2)- C_{12}\log (\|s\|^2(o) +\e)  \tag4.6
$$
for some constants $C_{11}$ and $C_{12}$ independent of $r$.
Since $s(o)\ne 0$   by the construction as one can specify the value of $s(o)$ and since the set $\{s=0\}$ is of measure zero, letting $\e\to 0$, the proof of the second statement in (ii) is completed.

If we only assume that $u$ is of at most linear growth, then using similar method, instead of (4.6), we have that
$$
\int_0^r \sigma\lf(\aint_{B_{o}(\sigma)}{\Cal R}(x)\cdot
\frac{\|s\|^2(x)}{\|s\|^2(x)+\e}\, dx \ri)\, d\sigma \le C \lf(r- {\Cal R}(o)\cdot
\frac{\|s\|^2(o)}{\|s\|^2(o)+\e} \ri)
$$
for some constant $C$ independent of $r$. The result follows by
letting $\e\to0$ as before.
\enddemo

Our first result on the structure of complete K\"ahler manifolds
with nonnegative bisectional curvature is a splitting theorem in
terms of harmonic function and holomorphic function.
Together with Lemmas 4.1 and 4.2, this theorem will be used from
time to time in the rest of this section.

\proclaim{Theorem 4.1} Let $M^m$ be a complete noncompact K\"ahler
manifold with nonnegative holomorphic bisectional curvature.
Suppose $f$ is a nonconstant harmonic function on $M$ such that
$$
\limsup_{x\to\infty}\frac{|f(x)|}{r^{1+\e}(x)}=0, \tag4.7
$$
for any $\e>0$, where $r(x)$ is the distance of $x$ from a fixed
point. Then  $f$ must be of linear growth and  $M$ splits
  isometrically as $\widetilde{M}\times\R$.
 Moreover the universal cover $\overline {M}$ of $M$  splits
  isometrically and holomorphically as $\widetilde{M'}\times\C$,
where $\widetilde M'$ is a complete K\"ahler
  manifold with nonnegative holomorphic bisectional curvature. Suppose that
there exists a nonconstant holomorphic function $f$ on $M$ satisfying (4.7). Then $M$ itself splits as $\widetilde{M}\times\C$.
\endproclaim

We need the following lemmas for the proof of Theorem 4.1.

The first one is a result in [L1, Corollary 5]. For the sake of
completeness, we will sketch the proof. It seems that in the proof
of this result, we need to assume that the holomorphic bisectional
curvature is nonnegative.

\proclaim{Lemma 4.3} Let $M$ be a complete noncompact K\"ahler
manifold with nonnegative holomorphic bisectional curvature. If
$f$ is a harmonic function with sub-quadratic growth defined on
$M$, then $f$ is pluri-harmonic.
\endproclaim
\demo{Proof} Let
$h=||f_{\abb}||^2=g^{\a\bar\delta}g^{\gamma\bar\beta}f_{\abb}f_{\gamma\bar\delta}$,
 where $g_{\abb}$ is the metric of $M$ and $g^{\abb}$ is its
inverse. Since $f$ is harmonic, by Lemma 2.1 in \cite{NT 2}satisfies (2.1) we know that
$$
\Delta
f_{\gamma\bar\delta}=-R_{\beta\bar\a\gamma\bar\delta}f_{\abb}+\frac12\lf(R_{\gamma
\bar p}f_{p\bar\delta}+R_{p\bar\delta}f_{\gamma\bar p}\ri).
$$
Hence in normal coordinates so that at a point $x$,
$f_{\abb}=\lambda_{\a}\delta_{\a\b}$, we have
$$
\split
\Delta h & = 2f_{\gamma\bar\delta s\bar s}f_{\delta\bar \gamma}+||f_{\abb\gamma}||^2+||f_{\a\bb\bg}||^2\\
&=-2R_{\b\a\gamma\bar\delta}f_{\abb}f_{\delta\bar\gamma}+\lf(R_{\gamma \bar p}f_{p\bar\delta}+R_{p\bar\delta}f_{\gamma\bar p}\ri)f_{\delta\bar\gamma}+||f_{\abb\gamma}||^2+||f_{\a\bb\bg}||^2\\
&= -2R_{\a\ba\gamma\bar\gamma}\lambda_{\a}\lambda_{\gamma}+2R_{\gamma\bg}\lambda_\gamma^2+||f_{\abb\gamma}||^2+||f_{\a\bb\bg}||^2\\
&=\sum_{\a,\b}R_{\a\ba\b\bb}\lf(\lambda_{\a}-\lambda_{\b}\ri)^2+||f_{\abb\gamma}||^2+||f_{\a\bb\bg}||^2
\\
&\ge 0,
\endsplit
$$
where we have used the fact that $M$ has nonnegative holomorphic
bisectional curvature. Since $|f(x)|=o\lf(r^2(x)\ri)$ where $r(x)$
is the distance from a fixed point $o\in M$, as in [L1, p.90-91],
we have
$$
\frac{1}{V_o(R)}\int_{B_o(R)}h\le
\frac{C}{R^{-2}V_o(R)}\int_{B_o(R)}|\nabla f|^2=o(1),
$$
as $R\to\infty$. Here   $C$ is a constant independent of $R$ and
we has used the gradient estimate in [C-Y]. Since $h$ is
subharmonic, $h\equiv0$ by the mean value inequality in [L-S].
Hence $f$ is pluri-harmonic.
\enddemo

\proclaim{Lemma 4.4} Let $M$ be a complete noncompact K\"ahler
manifold with nonnegative holomorphic bisectional curvature. Let
$f$ be a pluri-harmonic function. Then $\log (1+|\nabla f|^2)$ is
pluri-subharmonic.
 \endproclaim
\demo{Proof} We adapt the complex notation. Let $h=|\nabla
f|^2=\sum_{\a,\b}g^{\abb}f_{\a}f_{\bar{\b}}$. Here $g_{\abb}$ is
the K\"ahler metric and $(g^{\abb})$ is the inverse of
$(g_{\abb})$. To prove that $\log (1+h)$ is pluri-subharmonic, it
is sufficient to  show that $\lf[\log (1+h)\ri]_{\g\bg}\ge 0$ in
 normal coordinates. Direct calculation shows that:
$$
\split
h_{\g\bar{\g}} & =  \left(\sum_{\a\b}g^{\abb}f_{\a}f_{\bar{\b}}\right)_{\g\bar{\g}}\\
& =
\sum_{\a,\b}g^{\abb}\lf[f_{\a\g}f_{\bar{\b}\bar{\g}}+f_{\a\bar\g}f_{\bar\b\g}+f_{\a\g\bar{\g}}f_{\bar{\b}}
+f_{\a}f_{\bar{\b}\g \bar{\g}}\ri]\\
& = \sum_{\a}  f_{\a\g}f_{\bar{\a}\bar{\g}}+ \sum_{\a,s}
R_{\g\bar{\g}\a \bar{s}}f_{s}f_{\bar{\a}}
\endsplit \tag 4.8
$$
where we have used the fact that $f$ is pluri-harmonic. Hence
$$
\split \lf[\log (1+h)\ri]_{\g\bg}&=\frac{1}{(1+h)^2}
\lf[(1+h)h_{\g\bg}-h_\g h_\bg\ri]\\
&=\frac{1}{(1+h)^2}\bigg[(1+h)\lf(\sum_{\a}  f_{\a\g}f_{\bar{\a}\bar{\g}}+ \sum_{\a,s}  R_{\g\bar{\g}\a \bar{s}}f_{s}f_{\bar{\a}}\ri)\\
& \quad-\sum_\a  f_{\a\g}f_{\ba}\sum_\a f_\a f_{\ba\bg} \bigg]\\
&\ge\frac{1}{(1+h)^2} \lf(\sum_{\a}  f_{\a\g}f_{\bar{\a}\bar{\g}}+
\sum_{\a,s}  R_{\g\bar{\g}\a \bar{s}}f_{s}f_{\bar{\a}}\ri)
\endsplit
\tag 4.9
$$
where we have used the fact that $f$ is pluri-harmonic. From
(4.9), the fact that  $M$ has nonnegative
 holomorphic bisectional curvature, it is easy to see that $\log (1+h)$ is
pluri-subharmonic.
\enddemo

\demo{Proof of Theorem 4.1} Let $f$ be a nonconstant   harmonic
function on $M$ satisfying (4.7). Then $f$ is pluri-harmonic by
Lemma 4.3. By Lemma 4.4, the function $u=\log(1+|\nabla f|^2)$ is
pluri-subharmonic. By the gradient estimates in [C-Y],
$|u|(x)=o(\log r(x))$.
  By Theorem 3.2, we conclude that
 $|\nabla f|$ is constant. Hence $f$ must be of
linear growth. Moreover,  by the Bochner formula, we conclude that
$\nabla f$ must be parallel. Hence $J(\nabla f)$ is also parallel,
where $J$ is the complex structure of $M$.  From this it is easy
to see that the universal cover of $M$ splits as $\widetilde
M'\times\C$ isometrically and holomorphically. At the same time by
integrating along $\nabla f$, $M$ splits as $\widetilde M \times
\R$ isometrically, where $\widetilde M$ can be taken  to be the
component of $f^{-1}(0)$. In this case that $M$ supports a
nonconstant holomorphic function of growth rate (4.7), both the
real and imaginary part will split a factor of $\R$ and clearly
that they consist a complex plane $\C$.
\enddemo
  An easy consequence is that if the
  Ricci curvature is positive at some point of the manifold then any harmonic function satisfying (4.7) must
be a constant.

In the next theorem, we shall give some results on the  structure of
complete noncompact K\"ahler manifold with nonnegative sectional
or holomorphic bisectional curvature.

\proclaim{Theorem 4.2} Let $M^m$ be a complete noncompact K\"ahler
manifold with nonnegative holomorphic bisectional curvature.

\roster
\item"{(i)}" Suppose $M$ is simply connected, then $M=N\times M'$
holomorphically and isometrically, where $N$ is a compact simply connected K\"ahler manifold, $M'$ is a complete noncompact K\"ahler
manifold and both $N$ and $M'$ have  nonnegative holomorphic bisectional curvature.
Moreover, $M'$ supports a smooth strictly plurisubharmonic
function with bounded gradient and   satisfies {\rm ($\vgk$)} and {\rm ({\bf
CD})}, where $k$ is the complex dimension of $M'$. If, in addition,
$M$ has nonnegative sectional curvature outside a compact set,
then $M'$ is also Stein.

\item"{(ii)}"
If the holomorphic bisectional curvature of $M$ is positive at some point, then $M$ itself supports a smooth strictly plurisubharmonic function with bounded gradient, and satisfies {\rm($\text{\bf VG}_m $)} and {\rm({\bf CD})}, where $m$ is the complex dimension of $M$. If, in addition, $M$ has nonnegative sectional curvature outside a compact set, then $M$ is also Stein.
\endroster
\endproclaim
\proclaim{Remark 4.1}
\roster
\item"{(a)}" The factor $N$ in (i) may not be present. In
this case, $M=M'$ and satisfies the conditions on $M'$ mentioned in
the theorem. This kind of remark also applies to Theorem 4.3.
\item"{(b)}" It was first proved in \cite{CZ 2} that $M$ satisfies {\rm ($\text{\bf VG}_m $)} if the holomorphic bisectional curvature nonnegative and is positive at some point, and that $M$ satisfies {\rm({\bf CD})} if $M$ has positive holomorphic bisectional curvature everywhere.
\item"{(c)}" By \cite{M 2, HSW} (see also \cite{CC}), $N$ in (i) is a compact Hermitian symmetric manifold.
\endroster
\endproclaim
\demo{Proof} Let  $\Cal B$ be the Busemann function of $M$ and let $v$ be the solution of the heat equation with initial value $\Cal B$. Let $T_0$ be as in Lemma 4.1.

(i) Let $0<t<T_0$  be fixed and let $u(x)=v(x,t)$. By Lemma 4.1,
suppose $M$ is simply connected, then $M=N_1\times M_1$
isometrically and holomorphically so that $u_\abb\equiv0$ when
restricted on $N_1$ and $u_\abb>0$ when restricted to $M_1$.
Suppose $N_1$ is not compact,   then there is a ray of $M$ which
lies on $N_1$. By Lemma 4.1(iii), $u$ is not constant on $N_1$.
However, $u$ has bounded gradient by Lemma 4.1(ii). Theorem 4.1
then implies that $N_1=N_2\times \C$ isometrically and
holomorphically. Continuing in this way, we conclude that
$N_1=N\times \C^{\ell}$ isometrically holomorphically for some
$\ell\ge0$, where $N$ is a compact simply connected with
nonnegative holomorphic bisectional curvature.  Let
$M'=\C^\ell\times M_1$. Then $M'$ supports a strictly
plurisubharmonic function with bounded gradient and hence also
satisfies ($\vgk$) and ({\bf CD}) by Lemma 4.2, where $k=\dim_\C
M'$.

If, in addition, $M$ has nonnegative sectional curvature outside a
compact set, then $\Cal B$ is an exhaustion function by \cite{CG
2, GW 1}. Hence $u$ is an exhaustion function by Lemma 4.2.
Therefore $M'$ in the above is also Stein.

(ii) Suppose the holomorphic bisectional curvature of $M$ is
positive at some point, then $u$ is strictly plurisubharmonic by
Lemma 4.1(i). The rest of the proof is similar to (i).
\enddemo

In \cite{W 1}, Wu proved that a complete noncompact K\"ahler
manifold is Stein if  it has nonnegative sectional curvature
outside a compact set,    with nonnegative holomorphic bisectional
curvature everywhere which is positive outside a compact set. The
last statement of Theorem 4.2(ii) is a generalization of this
result.

In the last part of Theorem 4.1(i) or (ii), the assumption on the
sectional curvature is needed only for the proof that the Busemann
function is an exhaustion function. In some cases, this is true
even if we only assume that the Ricci curvature is nonnegative.
Hence we have the following result.

 \proclaim{Corollary 4.1} Let $M^m$ be a complete K\"ahler
manifold with nonnegative holomorphic bisectional curvature.
Suppose that   $M$ is of maximum volume growth or $M$ has a pole.   Then $M$ is Stein. Moreover, $M$ satisfies {\rm({\bf CD})}.
\endproclaim
\demo{Proof} Suppose $M$ has maximum volume growth. Let $\tilde M$ be the universal cover of $M$, then
$\tM$ has maximum volume growth and $\pi:\tM\to M$ is a finite
cover by \cite{L 1, p.10}. Suppose $\tM$ is Stein, then $\tM$ has a
smooth strictly plurisubharmonic exhaustion function $f$. Then the
function $h(x)=\sum_{\tx}f(\tx)$ for $x\in M$, where the summation
is taken over all $\tx\in \tM$ so that $\pi(\tx)=x$. Then $h$ is a
strictly  plurisubharmonic exhaustion function of $M$. Hence $M$
is also Stein. So without loss of generality, we may assume that
$M$ is simply connected.

Since $M$ has maximum volume growth, the Busemann function ${\Cal
B}$ is  an exhaustion function by \cite{Sh, p.400-401}. Let $u$ be as in the
proof of Theorem 4.2, then by this theorem, $M=N\times M'$ as
described in the theorem. Since $M$ has maximum volume growth, the
factor $N$ will not be present. Hence $M=\C^\ell\times M_1$ holomorphically and isometrically, so that   $u$ is strictly
plurisubharmonic on $M_1$. By Lemma 4.1, it is also an exhaustion
function on $M_1$. Therefore $M_1$ must be Stein by \cite{G} and so $M$ is also Stein. The last statement  follows from Lemma 4.1.

Suppose $M$ has a pole, then it is easy to see that the Busemann
function with respect to the pole is an exhaustion function. The
manifold is  diffeomorphic to $\R^{2m}$. One can   conclude that that the splitting given by Theorem 4.2 contains no   compact factor.  One can then proceed as above to conclude that $M$ is Stein.
  \enddemo

In [W 2], it was proved that $M$ is Stein under the assumption
that $M$ has a pole and nonnegative bisectional curvature which is
positive outside a compact subset of $M$. Our result answers
affirmatively the question raised in \cite{W 2, page 255} for the
nonnegative bisectional curvature case. Under the maximum volume
growth assumptions, if the holomorphic bisectional curvature is
actually positive everywhere then it is easy to see that it is
Stein by the results on smooth approximation of strictly
plurisubharmonic function in \cite{GW 4} and the result in
\cite{Sh} mentioned above. This was observed in \cite{WoZ}. Under
the maximum volume growth and the nonnegativity of the bisectional
curvature assumption together with the additional assumption that
the curvature decays like $r^{-1-a}$, the result was proved in
\cite{CZ 1}. This kind of results are related to a conjecture by
Greene-Wu \cite{GW 3} and Siu \cite{Si} that a complete noncompact
K\"ahler manifold with positive bisectional curvature is Stein.

Without assuming that $M$ is simply connected or the holomorphic bisectional curvature of $M$ is positive at some point, by applying Theorem 4.1 inductively, we immediately have:

\proclaim{Corollary 4.2} Let $M^m$ be a complete noncompact
K\"ahler manifold with nonnegative holomorphic bisectional
curvature. Then $M$ have the holomorphic-isometric splitting
$M^m=\C ^{k}\times M_2^{m-k}$. Here $M_2$ is a complete K\"ahler
manifold of nonnegative bisectional curvature with the property
that any holomorphic function on $M_2$ satisfying (4.7) must be a
constant.
\endproclaim

There is an open question whether the ring of polynomial growth
holomorphic functions on a complete noncompact K\"ahler manifold
with nonnegative curvature is finitely generated, see \cite{Y 4,
p. 391, Problem 63}.  This motivates us to study   the factor $M'$
in Theorem 4.2(i) or $M$ in Theorem 4.2(ii) in more details. We
have the following further splitting.

 \proclaim{Theorem 4.3} Let
$M^m$ be a complete noncompact K\"ahler manifold with nonnegative
holomorphic bisectional curvature. Assume that $M$ supports a
smooth strictly plurisubharmonic function $u$ on $M$ with bounded
gradient.

\roster
\item"{(i)}" If $M$ is simply connected, then $M=\C^\ell\times M_1\times
M_2 $ isometrically and holomorphically for some $\ell\ge0$, where
$M_1$ and $M_2$ are complete noncompact K\"ahler manifold with
nonnegative holomorphic bisectional curvature such that any
polynomial growth holomorphic function on $M$ is independent of
the factor $M_2$, and any linear growth holomorphic function is
independent of the factor $M_1$ and $M_2$. Moreover, $M_1$
supports a strictly plurisubharmonic function of logarithmic growth
and  satisfies {\rm({\bf FCD})} and {\rm($\text{\bf VG}_{a}$)}, for any
$a<k+1$,  where $k=\dim_\C M_1$.

\item"{(ii)}" Suppose the holomorphic bisectional curvature of $M$ is
positive at some point, then either $M$ has no nonconstant
polynomial growth holomorphic function or $M$ itself satisfies
{\rm({\bf FCD})} and {\rm($\text{\bf VG}_{a}$)}, for any $a<m+1$.
\endroster
\endproclaim
\demo{Proof} (i) We prove this part of the theorem by induction on
the dimension of $M^m$.

Suppose $M$ does not support any nontrivial polynomial growth
holomorphic function, then we simply take $M_2=M$ and the factors
$\C^\ell$ and $M_1$ are not present. Suppose there is a nontrivial
polynomial growth holomorphic function $f$ on $M$.  Let $w =\log
(1+|f|^2 )$. It is easy to see that $w$ is a plurisubharmonic
function so that
$$
0\le w(x)\le C_1\log (r(x)+2)\tag4.10
$$
We can solve the Cauchy problem (1.6) with initial data $w(x)$.
Denote $\tilde{w}(x,t)$ to be the solution.  Note that $\tilde
w(\cdot, t)$ is nonconstant because $w$ is nonconstant.  We can
apply Theorem 3.1 again to conclude that there exists $t>0$ and a
parallel distribution ${\Cal K}$ which corresponding to the null
space of $\tilde{w}_{\abb}(x,t)$. Suppose $\dim \Cal K=0$, then
$\tilde w(\cdot, t)$ is strictly plurisubharmonic with logarithmic growth by Corollary  1.4. Then $M=\C^\ell\times M_1$ by Corollary 4.1 so that every linear growth holomorphic function on $M$ is independent of $M_1$.    $M_2$ is not present in this case. Moreover, $M_1$
satisfies ({\bf FCD}) and ($\text{\bf VG}_{a}$) by Lemma 4.2.

Suppose $\dim\Cal K>0$, then $M=N_1\times N_2$, such that $\tilde
w_\abb(\cdot,t)\equiv0$ when restricted on $N_1$,  $\tilde
w_\abb(\cdot,t)>0$ when restricted on $N_2$. They are simply
connected, complete  K\"ahler manifolds with nonnegative
holomorphic bisectional curvature.    $\dim_\C N_1=\dim\Cal K>0$,
but $\dim_\C N_1<\dim_\C M$. Otherwise, $\tilde w(\cdot,t)$ is
harmonic on $M$ and is constant by (4.8) and \cite{CY}. Hence the
dimensions of $N_1$ and $N_2$ are both less than $m$. They are
also noncompact because $M$ supports a strictly plurisubharmonic
function. Hence $N_1$ and $N_2$ are K\"ahler manifolds satisfy the
same conditions satisfied by $M$.   By induction hypothesis
$N_1=\C^{\ell_1}\times N_{1,1}\times N_{1,2}$ and
$N_2=\C^{\ell_2}\times N_{2,1}\times N_{2,2}$ isometrically
holomorphically, such that for $j=1, \ 2,$ every polynomial growth
holomorphic function on $N_j$ is independent of the factor
$N_{j,2}$ and every linear growth holomorphic function is
independent of the factor $N_{j,1}\times N_{j,2}$. $N_{j,1}$
satisfies ({\bf FCD}) and ($\text{\bf VG}_{k_j}$) where
$k_j=\dim_\C N_{j,1}$. Let $M_1=N_{1,1}\times N_{2,1}$,
$M_2=N_{1,2}\times N_{2,2}$ and $\ell=\ell_1+\ell_2$. Then
$M=\C^\ell\times M_1\times M_2$ isometrically holomorphically.
Since every polynomial (respectively linear) growth holomorphic
function on $M$ is still a polynomial (respectively linear) growth
holomorphic function when restricted on $N_1$ and $N_2$, hence the
splitting satisfies all the required conditions if we can prove
that $M_1$ also satisfies the required volume growth and curvature
decay conditions.

The volume growth condition is satisfied by $M_1$ because of the
corresponding volume growth conditions are satisfied by $N_{1,1}$
and $N_{2,1}$. Moreover, for $r>0$, if $\Cal R_1$, $\Cal R'$ and
$\Cal R''$ are the scalar curvatures of $M_1$, $N_{1,1}$ and
$N_{2,1}$ respectively, and if $B_{o_1}(s)$, $B_{o'}(s)$ and
$B_{o''}(s)$ are the geodesic balls of $M_1$, $N_{1,1}$ and
$N_{2,1}$ respectively where $o_1=(o',o'')$, then
$$
\aint_{B_{o_1}(s)}\Cal R_1\le C\lf(\aint_{B_{o'}(s)}\Cal R'+\aint_{B_{o''}(s)}\Cal R''\ri)
$$
for some constant $C$ depending only on the dimensions of $M_1$,
$N_{1,1}$ and $N_{2,1}$. Since $N_{j,1}$ satisfies ({\bf FCD}) for
$j=1,\ 2$, $M_1$ also satisfies  ({\bf FCD}).

From the above proof, it is easy to see that part (i) is true if
$m=1$. The proof of Theorem 4.3(i) is then completed.

(ii) If the holomorphic bisectional curvature of $M$ is positive
at some point, suppose $M$ supports no nonconstant polynomial
growth holomorphic function, then we are done. Otherwise, let $f$
be a nontrivial polynomial growth holomorphic function. Construct
$w$ and $\tilde w$ as in the proof of (i), then there exists
$t>0$, $\tilde w(\cdot,t)$ must be strictly plurisubharmonic by
Theorem 3.1(iii). One can proceed as in the proof of (i).
\enddemo
\proclaim{Remark 4.2} By the Theorems 4.2 and 4.3, in order to
study polynomial growth holomorphic functions on a complete
noncompact K\"ahler $M^m$ with nonnegative holomorphic bisectional
curvature which is either simply connected or has positive
holomorphic bisectional curvature at some point, we may assume
that $M$ satisfies  the curvature decay condition {\rm({\bf FCD})} and
the volume growth condition {\rm($\text{\bf VG}_a$)} for any $a<m+1$. Note that under a
very mild condition on the bound of the  scalar curvature, if it
decays faster than {\rm({\bf FCD})}, then manifold must be flat. We
shall discuss this problem in a later section.
\endproclaim

As a simple consequence of Theorem 4.2 and 4.3 we have the
following  uniformization type result.

\proclaim{Corollary 4.3} Let $M$ be a complete, simply-connected,
K\"ahler manifold with nonnegative holomorphic bisectional
curvature. Suppose that the volume growth of $M$ satisfies $
V_o(r)=o(r^2)$. Then $M$ is biholomorphic to $N\times \C$, where
$N$ is biholomorphic to a compact Hermitian manifold. If $M$
supports nonconstant holomorphic functions of polynomial growth
the same result holds if $V_o(r)=O(r^{a})$, for some $a<3$.
\endproclaim

\input amstex
\documentstyle{amsppt}
\magnification=1200 \hsize=13.8cm \catcode`\@=11
\def\NoLogo{\let\logo@\empty}
\catcode`\@=\active \NoLogo
\def\vgk{\text{\bf VG}_k}

\def\heatt{\lf (\Delta-\frac{\p}{\p t}\ri)}
\def\heat{\lf(\frac{\p}{\p t}-\Delta\ri)}
\def\fp{f_+}
\def\fm{f_-}
\def \b {\beta}
\def\i{\sqrt{-1}}
\def\Ric{\text{Ric}}
\def\lf{\left}
\def\ri{\right}
\def\bbar{\bar \beta}
\def\a{\alpha}
\def\ol{\overline}
\def\g{\gamma}
\def\e{\epsilon}
\def\p{\partial}
\def\delbar{{\bar\delta}}
\def\ddbar{\partial\bar\partial}
\def\dbar{\bar\partial}

\def\C{\Bbb C}
\def\R{\Bbb R}
\def\tx{\tilde x}
\def\vp{\varphi}
\def\tv{\tilde v}

\def\hv{\hat v}

\def\dbar{\bar\partial}
\def\ba{{\bar\alpha}}
\def\bb{{\bar\beta}}
\def\abb{{\alpha\bar\beta}}

\def\i{\sqrt {-1}}

\def\tN{\widetilde N}
\def\tL{\widetilde L}

\def\tW{\widetilde W}
\def\tM{\widetilde M}
\def\tpi{\widetilde \pi}
\def\s{\Cal S}
\def\tS{\widetilde {\Cal S}}
\def\hpi{\widehat \pi}
\def\hM{\widehat M}
\def\hU{\widehat U}

\def \D {\Delta}
\def\tK{\widetilde K}
\def\pr{\text{\rm Proj}}
\def\aint{\frac{\ \ }{\ \ }{\hskip -0.4cm}\int}
\documentstyle{amsppt}
\vsize=19.0 cm

\subheading{\S5 Structure of non-negatively curved K\"ahler
manifolds II} \vskip.2cm

In \cite{Ta}, Takayama proved that if $M$ is a complete noncompact
K\"ahler manifold with nonnegative holomorphic bisectional
curvature and negative canonical bundle and if $M$ supports a
continuous plurisubharmonic exhaustion function, then $M$ has a
structure of holomorphic fibre bundle over a Stein manifold whose
fibre is biholomorphic to some compact Hermitian symmetric
manifold. In particular, the result applies to K\"ahler manifolds
with  nonnegative sectional curvature and positive Ricci
curvature.  This settled a conjecture of Greene-Wu \cite{GW 3,
\S3} that a complete noncompact K\"ahler manifold with nonnegative
sectional curvature and positive Ricci curvature is
holomorphically convex. In this section, we shall give   more
detailed results on the structure on the class of manifolds
related to the above conjecture. We shall also include  results of
Fangyang Zheng \cite{Z 2} on the structure of complete noncompact
K\"ahler manifold with nonnegative sectional curvature. The
authors are grateful to Fangyang Zheng for allowing them to
include his results and proofs in this work.

Before we state our results, let us make some preparations. Let
$M$ be a complete noncompact K\"ahler manifold  with nonnegative
holomorphic bisectional curvature.
 Let $\Cal B$ be the Busemann function at a point $o\in M$, and
 let $v$ be the solution to the heat equation with initial value
$\Cal B$. Then  there is $t_0>0$ such that $v(x,t_0)$ is smooth
plurisubharmonic and the kernel of $\Cal K(x,t_0)$ of
$v_\abb(x,t_0)$ is a smooth, parallel distribution on $M$.   In
the following, we shall suppress the variable $t_0$ and just write
$v(x)=v(x,t_0)$. Note that if $\Cal B$ is an exhaustion function of $M$, then $v(x)$ is also an exhaustion function. Moreover, $v$ has bounded gradient. All
these results are contained in  Lemma 4.1.

Let $\tM$ be the universal cover of $M$ with projection $\tpi$ and
let $\tv=v\circ \tpi$. Then $\tM=\tN\times \tL$ isometrically and
holomorphically. In the following, a point in $\tM$ will be
denoted by $(y,z)$ etc. The splitting of $\tM$ has the following
properties. For each $z\in \tL$, $v_\abb\equiv0$ on
$\tN_{z}=\tN\times\{z\}$ and for each $y\in \tN$, $v_\abb>0$ when
restricted on  $\tL_{y}=\{y\}\times \tL$. That is, $\tN_{z}$, for
$z\in \tL$ are integral manifolds of the distribution
$\lf(\tpi^{-1}\ri)_*(\Cal K)$.

Let $\Gamma$ be the fundamental group of $M$.

\proclaim{Lemma 5.1} Let $M$ be a complete noncompact K\"ahler manifold with nonnegative holomorphic bisectional curvature. Suppose that the Busemann function $\Cal B$ of $M$ is an exhaustion function and suppose the universal cover $\tM$ contains no Euclidean factor. Then with the above notations, $\tN$ is compact and a deck transformation $\gamma\in \Gamma$ is of the form $\gamma(y,z)=(\gamma_1(y),\gamma_2(z))$ so that $\gamma_1$ and $\gamma_2 $ are holomorphic isometries of $\tN$ and  $\tL$ respectively. Moreover $\gamma_2$ has no  fixed point unless $\gamma$ is the identity.
\endproclaim
\demo{Proof} Let us first  prove that $\tN$ is compact. Fix $z\in \tL$ and consider $\tN_z=\tN\times\{z\}$. By the construction of $\tN$, $\tv$ is pluri-harmonic on $\tN_z$. Suppose $\tv$ is not constant on $\tN_z$, then $\tN_z$ will contain a factor  $\C$ by Theorem 4.1. This contradicts the assumption that $\tM$ does not contain any Eucldiean factor. Hence $\tv$ must be constant on $\tN_z$. Since $\tv$ is the lift of $v$ and since $\Cal B$ and hence $v$ is an exhaustion function, we conclude that $\tpi(\tN_z)$ is a bounded and hence its closure $K$ in $M$ is compact.  Since $\tpi$ is a covering map,  there exists a compact set $\tK$ in $\tM$ such that $\tpi(\tK)\supset K$. We can now proceed as in \cite{CG 1, p.126}. Suppose $\tN_z$ is not compact, then there is a ray $\sigma$ in $\tN_z$. Since $\tpi(\sigma(n))\in K$, there exists $\gamma_n\in \Gamma$ such that $\gamma_n(\sigma(n))\in \tK$. Since $\tK$ is compact, passing to a subsequence if necessary, we may assume that $\gamma_n(\sigma(n))\to p\in \tM$ and $(\gamma_n)_*\lf(\sigma'(n)\ri)\to \vec  w\in T_p(\tM)$. Let $\tau$ be the geodesic with $\tau(0)=p$ and $\tau'(0)=\vec  w$, then it is easy to see that $\tau$ is a line. By \cite{CG 1}, $\tM$ has a factor of $\R$. This is a contradiction. Hence $\tN$ is compact.

Let $\gamma\in \Gamma$, then $\pr_2\lf(\gamma(\tN_z)\ri)$ is a compact subvariety in $\tL$ where $\pr_2$ is the projection onto $\tL$. Since $\tL$ supports a strictly pluri-subharmonic function, it must be a point. Hence $\gamma$ is of the form $\gamma(y,z)=(f(y,z),g(z))$. Since $\gamma$ is an isometry, $g$ will not increase length. This is also true for $g^{-1}$, and hence it is easy to see that $g$ is a local isometry and $f(y,z)=f(y)$.
Therefore $\gamma(y,z)=(\gamma_1(y),\gamma_2(z))$, where $\gamma_1$ and $\gamma_2$ are holomorphic isometries on $\tN$ and $\tL$ respectively. Suppose $\gamma$ is not the identity and suppose $\gamma_2(z)=z$ for some $z\in \tL$, then $\gamma_1$ will not have any fixed point. Then $\tN$ will cover a compact complex manifold with nonnegative holomorphic bisectional curvature. This is impossible by \cite{HSW}. This completes the proof of the lemma.
\enddemo
Let $\Gamma_2$ be the subgroup of the isometry group of $\tL$ which is the image under the map $\gamma\to \gamma_2$. By the lemma $\Gamma_2$ acts freely and holomorphically on $\tL$. Let $\hM=\tL/\Gamma_2$. Since $M=\lf(\tN\times \tL\ri)/\Gamma$, by the lemma there is a projection $\pi_r:M\to \hM$ such that the following diagram commutes:

$$\CD
\tN\times \tL @>\pr_2>> \tL\\
@V\tpi VV  @V\hpi VV \\
M @>\pi_r>> \hM.
\endCD
$$
In fact, $\pi_r(\tpi(y,z))=\hpi(z)$.
\proclaim{Theorem 5.1}   Let $M$ be a complete noncompact K\"ahler manifold   with nonnegative holomorphic bisectional curvature such that the Busemann function is an exhaustion function. Suppose the universal cover $\tM$ has no Euclidean factor. Using the above notations, we have the following:
\roster
\item"{(i)}" $\tM=\tN\times \tL$ and $\tN$ is compact.
\item"{(ii)}" $\pi_r:M\to\hM$ has the structure of a holomorphic fibre bundle, where each fibre is isometrically biholomorphic to $N$.
\item"{(iii)}" $\hpi:\tL\to \hM$ is a holomorphic and Riemannian covering map.   $\hM$ is a complete noncomapct K\"ahler and Stein manifold with nonnegative holomorphic bisectional curvature. $\tL$ is also Stein.
\endroster
\endproclaim
\demo{Proof} (i) follows from Lemma 5.1.

To prove (ii), let $z\in \tL$ and let $\hU$ be a neighborhood of $\hpi(z)$ in $\hM$ which is evenly covered by a family $\Cal F$ of neighborhoods  in $\tL$. Let $\tW_1$ and $\tW_2$ be two members in $\Cal F$, then there exists $\gamma_2\in \Gamma_2$ such that $\gamma_2(\tW_1)=\tW_2$. Suppose $\gamma_2$ is such that $\gamma=(\gamma_1,\gamma_2)$ for some $\gamma\in \Gamma$. Then $\gamma(\tN\times \tW_1)=\tN\times \tW_2$ and $\tpi$ maps $\tN\times \tW_1$ holomorphically and isometrically onto its image. Let $\tW_1$ be   one of the neighborhoods in $\Cal F$.   It is easy to see that $\pi_r^{-1}( \hU)= \tpi(\bigcup_{\tW\in\Cal F}\tN\times \tW)=\tpi(\tN\times \tW_1)$ which is holomorphically isometric to $\tN\times \tW_1$. Each fibre $\pi_r^{-1}(\hpi(z))=\tpi(\tN_z)$ which is holomorphically isometric to $\tN$. This completes the proof of (ii).

(iii) Let $\tv$ be the smooth pluri-subharmonic function on $\tM$ defined before Lemma 5.1. Since $\tN$ is compact, $\tv(y,z)=\tv(z)$ which is independent of the factor $\tN$. Since $\tv=v\circ \tpi$, for any $\gamma \in \Gamma$, $\tv(y,z)=\tv(\gamma(y,z))$. Hence $\tv$ is equivariant with respect to $\Gamma_2$ and so it descends to be a smooth strictly pluri-subharmonic function $\hv$ on $\hM$ because $\tv$ is strictly pluri-subharmonic on $\tL$. Note also that $v(x)=\hv(\pi_r(x))$. Since $v$ is an exhaustion function on $M$, $\hv$ is an exhaustion function of $\hM$. Hence $\hM$ is Stein by \cite{G}. The fact that $\tL$ is Stein follows from a result of \cite{St} that the universal cover of a Stein manifold is Stein.
\enddemo

\proclaim{Remark 5.1} The condition in the theorem that   $\Cal B$ is an exhaustion function will be satisfied if  $M$ has nonnegative sectional curvature outside a compact set, see \cite{CG 2, GW 1}. The condition that $\tM$ has no Euclidean factor will be satisfied if the Ricci curvature of $M$ is positive at some point. Hence if $M$ has nonnegative sectional curvature outside a compact set (in addition to the fact that $M$ has nonnegative holomorphic bisectional curvature) and if $M$  has positve Ricci curvature is   at some point, then $M$ has the structure as described in the theorem.
\endproclaim

Suppose $M$ has nonnegative sectional curvature everywhere, then Fangyang Zheng \cite{Z 2} obtains the  following stronger structure theorem.

\proclaim{Theorem 5.2}\cite{\rm Zheng} Let $M$ be a complete noncompact K\"ahler
manifold with nonnegative sectional curvature. Assume that
the universal cover $\tM$ of $M$ does not have Euclidean factors.
Then $M$ is simply connected and $M=N\times L$ isometrically and
holomorphically where $N$ is a compact Hermitian symmetric
manifold and $L$ is diffeomorphic to the Euclidean space
$\R^{2\ell}$ where $\ell=\dim_\C L$.
\endproclaim

The theorem follows immediate from the following lemma and Theorem 5.1.

\proclaim{Lemma 5.2} Let $M^m$ be a complete noncompact K\"ahler manifold with nonnegative sectional curvature so that its universal cover does not contain a Euclidean factor. Suppose that $M$ supports a stictly pluri-subharmonic function. Then the soul of $M$ is a point and hence $M$ is simply connected and is diffeomorphic to $\R^{2m}$.
\endproclaim

Let us assume the lemma is true of the moment.

\demo{Proof of Theorem 5.2} Apply Theorm 5.1 to $M$.  Using
the notations as in Theorem 5.1,  the manifold $\hM$ satisfies all
the conditions in the lemma because $\tL$ contains no Euclidean factor. Hence $\hM$ is simply connected. So $M$
is also simply connected because it is a fibre bundle over $\hM$ with
fibre $N$ which is simply connected. Hence $M=\tM$. By the lemma
again, we know that $\tL$ is diffeomorphic to the Euclidean space.
\enddemo

It remains to prove  Lemma 5.2.

\demo{Proof of Lemma 5.2} Let $\tM$ be the universal cover of $M$ with covering map $\pi$ and let $\s$ be a soul of $M$. Let $\tS=\pi^{-1}(\s)$.  Since $\s$ is totally geodesic and totally convex by \cite{CG 2}, $\tS$ is also totally geodesic and totally convex. In particular, $\tS$ is connected.  Suppose $\tS$ is noncompact, then it contains a ray which will also be a ray in $\tM$. Since $\pi(\tS)=\s$ which is compact, we can conclude as in the proof of Lemma 5.1 that $\tS$ is compact because $\tM$ contains no Euclidean factor. Morover, $\tS$ is simply connected by \cite{CG 2, Theorem 2.1}. We want to prove that $\tS$ is a point. It will be sufficient to prove that $\tS$ is flat because $\tS$ is connected, simply connected and  compact.

By \cite{Theorem 3.1, CG 2}, since $\s$ is a soul of $M$ and
$\tS=\pi^{-1}(\s)$, we have
$$
R(u,v)v=R(v,u)u=0\tag5.2
$$
for any point $p\in \tS$, any vector  $u\in T_p(\tS)$ and $v$ which is normal to $T_p(\tS)$. Here $R$ is the Riemannian curvature tensor of $\tM$. Let $J$ be the complex structure of $M$. Let $p\in \tS$ and let $W$ be the subspace of $V=T_p(\tS)$ consisting of  vectors $v\in V$ such that $Jv\in V$. Let $\gamma$ be a piecewise smooth closed curve on $\tS$ from $p$. Since $\tS$ is totally geodesic, $J$ commutes with parallel translation along $\gamma$ on $\tS$. Hence $W$ is invariant under parallel translation along $\gamma$. Note also that $J(W)=W$. Since $\tS$ is simply connected, $\tS$ can be splitted as a product with a factor $\tS_1$ whose tangent spaces are invariant under $J$.  By assumption, $M$ and hence $\tM$ supports a strictly pluri-subharmonic function. This implies that $\tS_1$ is a point and $W=\{0\}$. That is to say
$$
JV\cap V=\{0\}.\tag5.3
$$
For $v\in V$, let $A(v)$ be the orthorgonal projection of $Jv$
onto $V$.   Since $J+J^t=0$, it is easy to see that
$A+A^t=0$. We can find orthonormal basis $e_1,e_2,\dots,e_{2k-1},e_{2k},e_{2k+1},\dots,e_s $, $s=\dim_\R V$,   under which $A$ takes
the block diagonal form $\  A=\text{diag} \{ \delta_1J_2, \cdots ,
\delta_k J_2, 0_l\}$, where $k\geq 0$, $\ 1\geq \delta_1 \geq
\cdots \geq \delta_k >0$, $l=\dim (S) -2k$,
$$ J_2 = \left[
\matrix  0 & 1 \\ -1 & 0 \endmatrix \right],
 $$
and $0_l$ is a zero matrix. By (5.3), we have $\delta_j<1$.

Now take $e_1,\ e_2 $ for exmaple. We have
$$
Je_1=\a e_2+u_1, \ \ Je_2=-\a e_1+u_2,
$$
where $u_1,\ u_2\in V^\perp.$ Since $0<\a=\delta_1<1$,   $u_1$, $u_2$ are
nonzero vectors. Since $\langle Je_i,Je_j\rangle=\delta_{ij}$,
we conclude that    $\langle
u_i,u_j\rangle=\beta^2\delta_{ij}$, where
$\beta=||u_1||=||u_2||=\sqrt{1-\a^2}$. Let $e_3=\b^{-1}u_1$ and
let $e_4=\b^{-1} u_2$. Then we have

$$
Je_1=\a e_2+\b e_3, \ \ Je_2=-\a e_1+\b e_4.\tag 5.4
$$
 Using the fact that $J^2=-I$, from (5.4), we also have
$$
Je_3=-\b e_1-\a e_4,\ \ Je_4=-\b e_2+\a e_3.\tag 5.5
$$
 By (5.2) and (5.4), since $M$ is K\"ahler, we have
$$
0=-R(e_1,e_3,Je_4,e_1) =R(e_1,e_3, e_4,Je_1) =R(e_1,e_3, e_4, \a
e_2+\b e_3) =\a R(e_1,e_3,e_4,e_2).
$$
  By (5.2) and (5.4),
$$
0=-R(e_1,e_4,Je_2,e_1)=R(e_1,e_4, e_2,Je_1) =R(e_1,e_4, e_2, \a
e_2+\b e_3) =\b R(e_1,e_4,e_2,e_3).
$$
Hence by the Bianchi identity, we have
$$
R(e_1,e_2,e_3,e_4)=- R(e_1,e_3,e_4,e_2)-R(e_1,e_4,e_2,e_3)=0.
$$
So
$$
\split
0&=R (e_1,e_2,e_3,e_4)\\
&=R(e_1,e_2,Je_3,Je_4)\\
&=R(e_1,e_2,-\b e_1-\a e_4,-\b e_2+\a e_3)\\
&=R(e_1,e_2,-\b e_1 ,-\b e_2 )+R(e_1,e_2,-\b e_1 , \a e_3)+R(e_1,e_2,-\a e_4,-\b e_2 )\\
&\quad+R(e_1,e_2,-\a e_4 , \a e_3)\\
&=\b^2R(e_1,e_2,e_1,e_2)
\endsplit\tag 5.5$$
where we have used (5.2)and (5.5).   Using (5.2) and (5.4),
we have
$$
R(e_1,Je_1,e_1,Je_1)=R(e_1,\a e_2+\b e_3,e_1,\a e_2+\b
e_3)=\a^2R(e_1,e_2,e_1,e_2)=0.\tag5.6
$$
Now we use the following fact in \cite{Z 1}. Suppose $X$ is such that the
holomorphic sectional curvature $R(X,\bar{X}, X,\bar{X})=0$ and if the sectional curvature is also nonnegative, then  $K(X,\bar{X}, Y,
\bar{Y})=0$ for any $Y$. Hence (5.6) implies that the sectional curvature  $K(e_1, u)$ of the plane spanned by $e_1$ and   any tangent vector $u\in T_p(\tS)$  is zero. Similarly, we can prove that $K(e_j,u)=0$ for $1\le j\le 2k$. Since $Je_j\in V^\perp$ for $2k+1\le j\le s$,   $K(e_j,Je_j)=0$ by (5.2). Hence we have $K(e_j,u)=0$ for all $j$ and for all $u\in T_p(\tS)$. Since $p$ is any point on $\tS$, $\tS$ is flat. This completes the proof of the lemma.
\enddemo

In case the universal cover of $M$ contains some Euclidean factors, then we have the following result which is also due to Fangyang Zheng \cite{Z 2}.

\proclaim{Corollary 5.1} Let $M$ be a complete K\"ahler manifold
with nonnegative sectional curvature. Then its universal cover is of the form $\tM=\C^k\times \tN\times \tL $ where $\tN$ is a compact Hermitian symmetric manifold, $\tL$ is Stein and  $\tL$ contains no Euclidean factor.  Moreover, there exists a discrete
subgroup $\Gamma \subseteq I_h({\C}^k)$ which acts freely on
${\C}^k$, and group homomorphisms $\rho : \Gamma \rightarrow
I_h(\tN)$, $\tau : \Gamma \rightarrow I_h(\tL)$, such that $M$ is
holomorphically isometric to the quotient of $\tM$ by $\Gamma $
which acts on $\tM$   as
$$ \gamma (x,y,z) = (\gamma (x), \rho (\gamma )(y), \tau
(\gamma )(z)) $$ for any $\gamma \in \Gamma$. In particular, $M$
is a holomorphic and Riemannian fiber bundle with fiber $\tN\times
\tL$ over the flat K\"ahler manifold ${\C}^k/\Gamma $. Here $I_h(X)$ denotes the group of isometric biholomorphisms of a K\"ahler manifold $X$.
\endproclaim
\demo{Proof} By Theorems 4.2 and 4.3, it is easy to see that $\tM$ is of the form as described in the corollary. Note that $\tN$ or $\C^k$ may reduce to a point. Let $G$ be the fundamental group of $M$. Let $\beta\in G$, we claim that  $$\beta(x,y,z)= (f(x), g(y), h(z))$$
for  $(x,y,z)\in \C^k\times \tN\times \tL$.  denote the point of $\tM$ according to the
splitting. A priori $f=f(x,y,z)$. So are $g$ and $h$. As in the proof of Lemma 5.1,   we have $f=f(x,z)$ and $h=h(x,z)$ are
independent of $y$ since $\tN$ is a compact Hermitian symmetric manifold. The next observation  that     $h$ is also independent of $x$. Otherwise, there
exists $x_1$, $x_2$ in $\C^k$ and $z\in \tL$ such that $h(x_1, z)  \ne h(x_2,z)$. Denote the line passing $x_1$ and $x_2$ in $\C^k$ to be
$\alpha(s)$. Then $\beta(\alpha)$ is also a line. This in
particular implies that $h(\a)$ is also a geodesic and
distance realizing, therefore a line, since $h(\a)$ is not a
point. This contradicts the assumption that $\tL$ does not contain
any lines. Hence $h=h(z)$. As in the proof of Lemma 5.1, we conclude that $f=f(x)$. Moreover,   it is easy to see that $f\in I_h(\C_k)$, $g\in I_h(\tN)$ and $h\in I_h(\tL)$.

Let $\rho_1:G\to I_h(\C^k)$ be the homomorphism defined by the above correspondence $\beta\to f$. Define $\rho_2$, $\rho_3$ similarly. We claim that $\rho_1$ is a monomorphism.
Otherwise, we can find $\beta\neq$identity in $G$ such that
 $f=$identity, where $\beta=(f,g,h)$. Then $(f,g)$ will act freely on $\tN\times \tL$. This
implies that the group generated by $(f,g)$ will be the fundamental group of some complete
K\"ahler manifold which is covered by $\tN\times \tL$ and is noncompact by \cite{HSW} or \cite{CG 1} because $\tN\times\tL$ contains no Euclidean factors. This contradicts
Theorem 5.2. Therefore we know that $\rho_1$ is an isomorphism.

Now simply denote $\Gamma =\rho_1(\tpi)$. Let $\rho=\rho_2\circ
\rho_1^{-1}$ and $\tau=\rho_3\circ\rho_1^{-1}$. This completes the proof of the corollary.
\enddemo

\input amstex
\documentstyle{amsppt}
\magnification=1200 \hsize=13.8cm \catcode`\@=11
\def\NoLogo{\let\logo@\empty}
\catcode`\@=\active \NoLogo

\def\heatt{\lf (\Delta-\frac{\p}{\p t}\ri)}
\def\heat{\lf(\frac{\p}{\p t}-\Delta\ri)}
\def\fp{f_+}
\def\fm{f_-}
\def \b {\beta}
\def\i{\sqrt{-1}}
\def\Ric{\text{Ric}}
\def\lf{\left}
\def\ri{\right}
\def\bbar{\bar \beta}
\def\a{\alpha}
\def\ol{\overline}
\def\g{\gamma}
\def\e{\epsilon}
\def\p{\partial}
\def\delbar{{\bar\delta}}
\def\ddbar{\partial\bar\partial}
\def\dbar{\bar\partial}

\def\C{\Bbb C}
\def\R{\Bbb R}
\def\tx{\tilde x}
\def\vp{\varphi}

\def\tN{\tilde\nabla}

\def\dbar{\bar\partial}
\def\ba{{\bar\alpha}}
\def\bb{{\bar\beta}}
\def\abb{{\alpha\bar\beta}}

\def\i{\sqrt {-1}}

\def \D {\Delta}
\def\aint{\frac{\ \ }{\ \ }{\hskip -0.4cm}\int}
\documentstyle{amsppt}
\vsize=19.0 cm

\subheading{\S6 Poincar\'e-Lelong equation and a gap theorem}\vskip .2cm

In this section, we shall solve the Poincar\'e-Lelong equation using some refined estimates developed in previous sections, in particular in \S1--\S3. One of the   motivation is  to  discuss  the  curvature decay condition ({\bf FCD}) defined in section 4, see Remark 4.2. We shall prove the following:

\proclaim{Theorem 6.1} Let $M^m$ be a complete noncompact K\"ahler manifold
with nonnegative holomorphic bisectional curvature.  Let $\rho$ be
a real closed (1,1)  form with trace $f$. Suppose $f\ge0$ and
$\rho$ satisfies  the following conditions:

$$
 \int_0^\infty  \aint_{B_o(s)}||\rho||ds<\infty,
\tag 6.1
$$
 and
$$
\liminf_{r\to\infty}\lf[\exp\lf(-ar^2\ri)\int_{B_o(r)}||\rho||^2\ri]<\infty
\tag 6.2
$$
for some $a>0$.
 Then there is a solution $u$ of
the Poincar\'e-Lelong equation $\sqrt{-1}\ddbar u =\rho$.
Moreover, for any $0<\epsilon<1$, $u$ satisfies
$$
\split
 \alpha_1r\int_{2 r}^\infty k(s)ds+\beta_1\int_0^{2r}sk(s)ds
&\ge u(x)\\
&\ge-\alpha_2r\int_{2r}^\infty k(s)ds-\beta_2\int_0^{\epsilon
r}sk(x,s)ds+\beta_3\int_{0}^{2r}sk(s)ds
\endsplit \tag 6.3
$$
for some positive constants $\alpha_1(m)$, $\alpha_2(m,\epsilon)$
and $\beta_i(m)$, $1\le i\le 3$, where $r=r(x)$. Here  $k(x,s)= \aint_{B_x(s)}f$ and $k(s)=k(o,s)$, where $o\in M$
is a fixed point.
\endproclaim

The theorem  was first proved in \cite{MSY} under the assumption that $M$ has maximal volume growth and $||\rho||(x)$ decays like $r^{-2}(x)$ pointwisely. Later in  \cite{NST 1, Theorem 5.1} the theorem was generalized  by assuming the following   the condition instead of (6.2):
$$
\liminf_{r\to\infty}\aint_{B_o(r)}||\rho||^2=0. \tag6.4
$$
(6.4) is obviously much stronger than (6.2). However, it would be nice if (6.2) can be totally removed.

We shall use the ideas in \cite{MSY} and [NST 1]. By [NST 1,
Theorem 1.1 and Theorem 1.2], there exists a solution to the
Poisson equation $\D u=f$ such that $u$ satisfies (6.3). The
difficult part is to prove that $\i\ddbar u=\rho$. As in
\cite{NST 1 }, one can prove that $||\i\ddbar u||$ behaves like
$||\rho||$ in the average sense. If (6.4) is assumed, the result
will follow easily by using the mean value theorem of \cite{LS, p. 287}.
If we only assume (6.2), the method does not work
because the average of $||\rho||$ might grow exponentially.

The outline of the proof of $\i\ddbar u=\rho$ is as
follows. First we solve  the Cauchy problem (1.6) with initial
data $u(x)$ for all time and let  $v(x,t)$ be the solution. Let
$w=u-v$. By an argument as in Lemma 2.1, one can show that
$\|\rho-\sqrt{-1}\ddbar w\|$ is a subsolution of the heat equation,
and that  $\|\rho-\sqrt{-1}\ddbar w\|(x,t)\to 0$ as $t\to \infty$
using (6.1) and (6.2). On the other hand, we shall prove that
$v(x,t)-v(o,t)$ together with its derivatives uniformly converges
to a constant over any fixed compact subset, at least
subsequentially. Therefore $\|\sqrt{-1}\ddbar v\|(x,t)\to 0$,
which implies that $\rho-\sqrt{-1}\ddbar u\equiv 0$.

As in \cite{NST 1}, by taking $M\times \C^2$, we may assume that
$M$ is nonparabolic and its Green's function $G(x,y)$ satisfies
(1.16), with $n=2m$ being the real dimension of $M$.

 As mentioned above, by \cite{NST 1}
 we can solve $\Delta u=f$ with $u$ satisfying (6.3). $u$ is given by
$$
u(x)=\int_M\lf(G(o,y)-G(x,y)\ri)f(y)dy. \tag6.5
$$
The details of the proof that $\i\ddbar u=\rho$ are contained in the following two lemmas.

\proclaim{Lemma 6.1} \roster
\item"{(i)}" The Cauchy problem (1.6) with initial
value $u$ has long time solution $v(x,t)$ which is given by
$$
v(x,t)=\int_MH(x,y,t)u(y)dy.
$$

\item"{(ii)}" There exists $t_i\to\infty$ such that $v(x,t_i)-v(o,t_i)$
together with their derivatives converge uniformly on compact
subsets to a constant function.
\endroster
\endproclaim
\demo{Proof} (i) We want to apply Lemma 1.2. For any $R>0$ and $x\in B_o(R)$,
$$
\split
|u(x)|&\le\int_M\lf|G(o,y)-G(x,y)\ri|f(y)dy\\
&=\lf\{\int_{M\setminus B_o(4R)}+\int_{B_o(4R)} \ri\}\lf|G(o,y)-G(x,y)\ri|f(y)dy\\
&=I(x)+II(x).
\endsplit\tag 6.6
$$
By (6.1), we have
$$
I(x)\le C_1r(x)\int_{2R}^\infty k(s)ds\le C_2r(x) \tag 6.7
$$
as in [NST 1, (1.4)], for some constants $C_1$ and $C_2$
independent of $R$. On the other hand, as in  \cite{NST 1, p. 347
and p. 356}, we have
$$
\split
\int_{B_o(R)}II(x)\, dx &\le \int_{x\in B_o(r)}\lf[\int_{y\in B_o(4R)}\lf(G(o,y)+G(x,y)\ri)f(y)dy\ri]dx\\
&=\int_{y\in B_o(4R)}\lf[\int_{x\in B_o(R)}\lf(G(o,y)+G(x,y)\ri)dx\ri]f(y)dy\\
&=C_3V_o(R)\lf(R^2k(4R)+\int_0^{4R}s k(s)ds\ri)
\endsplit
$$
 where $C_3$ depends only on $n$. Combining this with (6.6) and (6.7) and using (6.1), we  conclude that
$$
\aint_{B_o(R)}|u|\le C(1+R^2)
$$
for some constant independent of $R$. By Lemma 1.2, (i) follows.

(ii) Let us first give an estimate of $|\nabla v|$. We cannot use
the same method as in Lemma 1.5, because  we  have neither the
estimate of the integral of $u^2$ as in Lemma 1.5, nor the uniform
bound on $|\nabla u|^2$. However, we may proceed as in the proof
of Lemma 1.2. Namely, use cutoff functions $\vp_i$ and denote
$f_i=\vp_i f$. Solve $\D u_i=f_i$ by using (6.5) and find solution
$v_i$ of (1.6) with initial value $u_i$. Then $v_i$ subconverge to
$v$ together with their derivatives uniformly on compact sets of
$M\times[0,\infty)$. Note that $|\nabla u_i|$ is bounded by
\cite{NST 1, Theorem 1.3} and hence $|\nabla v_i|$ is bounded by
Lemma 1.5 or \cite{LT}. We can apply the maximum principle to
$|\nabla v_i|$ which is a subsolution of the heat equation and conclude that for any
$x$ such that $r(x)\le \sqrt t$,
$$
\split
|\nabla v_j|(x,t)&\le \int_M H(x,y, t)|\nabla u_i|(y)\, dy\\
&\le C_4 \sup_{r\ge
\sqrt{t}}\aint_{B_x(r)}|\nabla u_i|(y)\, dy\\
&\le C_4 \sup_{r\ge \sqrt{t}}\aint_{B_o(2r)}|\nabla u_i|(y)\, dy\\
&\le C_5\sup_{r\ge
\sqrt{t}}\lf(\int_{4r}^\infty k(s)\, ds +rk(4r)\ri)\\
& \le C_6 \int_{4\sqrt{t}}^\infty k(s)\, ds
\endsplit \tag6.8
$$
for some constants $C_4-C_6$ depending only on $n$. Here we have
used   Corollary 3.2 of [N 2]  in the second inequality,   Theorem
1.3 of [NST 1] in the fourth inequality and we have also used the
volume comparison as well as the fact that $0\le f_i\le f$. Hence
$$
\sup_{x\in B_o(\sqrt t)}|\nabla v_i|(x,t)\le C_6\int_{4\sqrt{t}}^\infty k(s)\, ds.
$$
for all $i$ and so
$$
\sup_{x\in B_o(\sqrt t)}|\nabla v|(x,t)\le C_6\int_{4\sqrt{t}}^\infty k(s)\, ds.\tag6.9
$$

On the other hand, $f_i$ has compact support, $u_i$ and $v_i$
are bounded. Since $ (v_i)_t$ is a   solution to the heat equation
with initial value $f_i$, as in the proof of Lemma 1.5 (or (6.12)
below), one can prove that for any $T>0$, there exist   constants
$C_i$ such that
$$
\int_0^T\aint_{B_o(r)}|(v_i)_t|^2\le C_i
$$
for all $r$. Hence we can apply maximum principle and conclude that
$$
\split
\frac{\p v_i}{\p t}(x,t)&=\int_M H(x,y,t) f_i(y)\, dy\\
&\le C(n)\sup_{r\ge \sqrt{t}}\aint_{B_x(r)} f_i(y)\, dy\\
&\le C(n)\sup_{r\ge \sqrt{t}} k(x,r).
\endsplit
$$
Here we also used Corollary 3.2 of [N 2]. Note we also have
$(v_i)_t\ge0$. Hence we have
$$
0\le \frac{\p v}{\p t}(x,t)\le C(n)\sup_{r\ge\sqrt{t}}k(x,r)\le
\frac{C(n)}{\sqrt{t}}\int_{\sqrt{t}}^\infty k(x,s)\, ds.\tag6.10
$$
By (6.9), (6.10) and (6.1), for any $t_0>1$, and $r>0$, the function $v(x,t)-v(o,t_0)$ is   bounded in $B_o(r)\times[t_0-1,t_0+1]$ by a constant which is independent of $t_0$ and  $\lim_{t\to\infty}\sup_{B_o(r)}|\nabla v(\cdot,t)|\to0$. Hence, it is easy to see that (ii) is true.
\enddemo
Now let $w=u-v$.
\proclaim{Lemma 6.3} As $t\to \infty$, $||\rho-\i\ddbar w||(x,t)$ converges to zero uniformly on compact subsets in $M$.
\endproclaim
\demo{Proof} We claim that
$$
\|\rho-\sqrt{-1}\ddbar w\|(x,t)\le \int_M H(x,y, t) \|\rho\|(y)\,
dy.\tag 6.11
$$
If (6.11) is true, then one can apply  Corollary 3.2 of [N 2] again to conclude that, for $x\in
B_o(\sqrt{t})$,
$$
\split \|\rho-\sqrt{-1}\ddbar w\|(x,t)& \le C(n)\sup_{r\ge
\sqrt{t}}\aint_{B_x(r)}\|\rho\|(y)\, dy\\
&\le C(n)\sup_{r\ge \sqrt{t}}\aint_{B_o(2r)}\|\rho\|(y)\,
dy \endsplit
$$
for some constant $C(n)$ depending only on $n$.
From the assumption (6.1), this implies that
$\sup_{B_o(\sqrt{t})}\|\rho-\sqrt{-1}\ddbar w\|(\cdot,t)\to 0$ as
$t\to \infty$ and the lemma follows.

To prove (6.11), we first observe that since $\rho$ is real
$d$-closed $(1,1)$ form, locally it can be written as
$\sqrt{-1}\ddbar \Psi$.  Using $\Delta \Psi =f$, it is easy to see
that $\Psi-w$ satisfies the heat equation. Hence
$\eta=\rho-\sqrt{-1}\ddbar w=\sqrt{-1}\ddbar(\Psi-w)$  satisfies
(2.1) by Lemma 2.1 in \cite{NT 2}. (6.11) will follow from Lemma
2.2 provided $\eta$ satisfies (2.2) and (2.3).  By (6.1), since
$w\equiv0$ at $t=0$, it is easy to see that  $\eta$ satisfies
(2.2) for any $a>0$. Next, we estimate $|\nabla^2 v|^2$. Again, we
may obtain the estimates for $v_i$ first and let $i\to\infty$.
Hence, as in the proof of Lemma 3.1, for any $T>1$ and $r^2\ge T$,
using the first inequality in (6.8) one can prove that

$$
\split
\int_0^{T}\aint_{B_o(r)}|\nabla^2 v|^2&\le C_1\lf[\frac1{r^2}\int_0^T\aint_{B_o(2r)}|\nabla v|^2
+\aint_{B_o(2r)}|\nabla u|^2\ri]\\
&\le C_2\bigg[(T+1) \aint_{B_o(8r)}|\nabla u|^2(x)dx\\
&\qquad
+\int_0^Tt^{-2}\lf(\int_{8r}^\infty \exp\lf(-\frac{s^2}{20 t}\ri)s\aint_{B_o(s)}|\nabla u|(y)\, dy
 ds\ri)^2dt \bigg] \\
&\le C_3\bigg[(T+1)\aint_{B_o(8r)} |\nabla u|^2(x)dx\\
&\qquad+\int_0^T
\lf(\int_{4r}^\infty \exp\lf(-\frac{s^2}{20 t}\ri) d\lf(\frac{s^2}{t}\ri)\ri)^2 \bigg]dt\\
&\le C_4(T+1) \lf[\aint_{B_o(4r)}|\nabla u|^2(x)dx+1\ri]\\
&\le C_5(T+1)\lf[\lf(\int_{16 r}^\infty k(s)ds\ri)^2+r^2\aint_{B_o(16r)}||\rho||^2+1\ri]\\
&\le C_6(T+1)\lf[r^2\aint_{B_o(8r)}||\rho||^2+1\ri]
\endsplit\tag 6.12
$$
for some constants $C_1-C_6$ independent of $r$ and $T$. Here we
have used  Lemma 1.1 in the second inequality, (6.1) and Theorem
1.3 of \cite{NST 1} in third inequality, Theorem 1.3 of [NST 1] in
the fifth inequality. By Theorem 1.3 in \cite{NST 1} again, we
have
$$
\aint_{B_o(r)}|\nabla^2 u|^2\le C_7\lf[r^2\aint_{B_o(2r)}||\rho||^2+1\ri]
$$
for some constant $C_7$ depending only on $n$. Combining this with
(6.12), we conclude that
$$
\int_0^{T}\aint_{B_o(r)}||\rho-\i\ddbar w||^2 \le C_7(T+1)\lf[r^2\aint_{B_o(2r)}||\rho||^2+1\ri]
$$
for some constant $C$ independent of $T$ and $r$. Hence  by (6.2),
 $\eta$ also satisfies (2.3), and we can apply Lemma 2.2
to conclude that (6.11) is true.
\enddemo
Now Theorem 6.1 follows from Lemmas 6.1 and 6.2.

\medskip

Using Theorems 6.1 and 3.2, we can prove that under a mild condition, the scalar curvature of a nonflat complete noncompact
K\"ahler manifold with nonnegative holomorphic bisectional
curvature cannot  decay faster than ({\bf FCD}).

\proclaim{Corollary 6.1} Let $M$ be a complete noncompact K\"ahler manifold with nonnegative holomorphic bisectional curvature and let $\rho\ge 0$ be a $d-$closed real
$(1,1)$ form. Assume that $\rho$
satisfies (6.2). Then $\rho\equiv 0$, if
$$
\int_0^rs\lf(\aint_{B_o(s)}\|\rho\|(y)\, dy\ri)\, ds=o(\log r). \tag 6.13
$$
In particular, if the scalar curvature $\Cal R$ of $M$ satisfies (6.2) and (6.13)  with $||\rho||$ replaced by $\Cal R$, then $M$ must be flat.
\endproclaim
\demo{Proof} By Theorem 6.1, one can solve $\sqrt{-1}\ddbar
u=\rho$ where $u$ satisfies (6.3).     (6.13) then implies that $u(x)=o(\log r)$ and $u$ must be constant by Theorem 3.2. Hence $\rho\equiv0$. The last result follows from this by considering $\rho$ to be the Ricci form of $M$.

\enddemo

\proclaim{Remark 6.1} The gap theorem in the last part of the
corollary was first obtained in \cite{MSY} under the assumptions
that $M$ has maximum volume growth with curvature decays like
$r^{-2-\e}$ pointwisely. These implies   (6.4) and (6.13) are true uniformly for all $o\in M$ (with $||\rho||$ replaced by $\Cal R$). Later, using K\"ahler-Ricci flow of \cite{S}, it was generalized
in \cite{CZ 1} by only assuming  that (6.13) holds uniformly
for all $o\in M$. In order to use the K\"ahler Ricci flow, it was also assumed that   $\Cal R$ is bounded in \cite{CZ 1}. It is easy to see that if $\Cal R$ is bounded, then (6.13) will imply (6.4).  In   Corollary 6.1, $\Cal R$ might   grow exponentially. Moreover, we only need to assume that (6.13) holds at one point.
\endproclaim

In Theorem 4.3(ii), it is proved that if $M$ has nonnegative holomorphic bisectional curvature which is positive at some point and if $M$ supports a nonconstant  polynomial growth holomorphic function, then it satisfies ({\bf FCD}). The following result is a partial converse of this.

\proclaim{Corollary 6.2} Let $M$ be a complete noncompact K\"ahler manifold with nonnegative holomorphic bisectional curvature. Suppose $M$ satisfies the condition
({\bf FCD}) and suppose the  Ricci curvature of $M$ is
positive at a point. If the scalar curvature ${\Cal R}$ also satisfies
(6.2), then $M$ supports nonconstant holomorphic functions of polynomial
growth. Moreover, if $Ric(o)>0$, there exists $\{f_1,\cdots,
f_m\}$, holomorphic functions of polynomial growth such that they
form a local coordinates near $o$. In particular, there exists a
positive constant $\delta$ independent of $d$ such that
$$
\dim ({\Cal O}_{d}(M))\ge \delta d^{m}
$$
for $d$ sufficient large. Here ${\Cal O}_d(M)$ is the vector space consisting all
holomorphic functions $f$ such that $|f(x)|\le C(r(x)+1)^d$ for
some constant $C(f)$.
\endproclaim
\demo{Proof} By Theorem 6.1, we can solve the Poincar\'e-Lelong
equation $\sqrt{-1}\ddbar u=\Ric$, since ({\bf FCD}) implies
$\|\Ric\|$ satisfies (6.1). Moreover, by (6.3) we know that $u(x)$
satisfies $u(x)\le C\log (r(x)+2)$ for some $C$. The corollary then follows from rather standard argument, see \cite{M 1, N 1} for example.  In fact, let $\{z_1,
\cdots,z_m\}$ be the local coordinate near $o$. Let $h_i=\vp(x)
z_i$, where $\vp(x)$ is a cut-off function which has support
inside the local coordinate neighborhood. Let $\theta_i=\dbar
h_i$. Now apply Theorem 3.2 in \cite{N 1}, with $E$ being the
anti-canonical line bundle. We then have functions $\eta_i$ such
that $\dbar \eta_i=\theta_i$ and $\eta_i(o)=0$. Moreover the
$\eta_i$ satisfies the following estimate:
$$
\int_M |\eta_i|^2\exp(-Cu(x))<\infty. \tag 6.14
$$
It is easy to see that $f_i=\theta_i-\eta_i$ will be  holomorphic
functions such that $f_i=z_i$ near $o$. Moreover $f_i$ satisfies
(6.14). Applying the mean value inequality of [LS, p.287] as in the proof of Lemma
4.2 we conclude that $f_i$ are of polynomial growth. The second
claim of the corollary follows from simple dimension counting.
\enddemo

\input amstex
\documentstyle{amsppt}
\magnification=1200 \hsize=13.8cm \catcode`\@=11
\def\NoLogo{\let\logo@\empty}
\catcode`\@=\active \NoLogo
\def\tM{\tilde M}
\def\vgk{\text{\bf VG}_k}
\def\heatt{\lf (\Delta-\frac{\p}{\p t}\ri)}
\def\heat{\lf(\frac{\p}{\p t}-\Delta\ri)}
\def\fp{f_+}
\def\fm{f_-}
\def \b {\beta}
\def\i{\sqrt{-1}}
\def\Ric{\text{Ric}}
\def\lf{\left}
\def\ri{\right}
\def\bbar{\bar \beta}
\def\a{\alpha}
\def\ol{\overline}
\def\g{\gamma}
\def\e{\epsilon}
\def\p{\partial}
\def\delbar{{\bar\delta}}
\def\ddbar{\partial\bar\partial}
\def\dbar{\bar\partial}

\def\C{\Bbb C}
\def\R{\Bbb R}
\def\tx{\tilde x}
\def\vp{\varphi}

\def\tN{\tilde\nabla}

\def\dbar{\bar\partial}
\def\ba{{\bar\alpha}}
\def\bb{{\bar\beta}}
\def\abb{{\alpha\bar\beta}}

\def\i{\sqrt {-1}}

\def \D {\Delta}
\def\aint{\frac{\ \ }{\ \ }{\hskip -0.4cm}\int}
\documentstyle{amsppt}
\vsize=19.0 cm

\enddocument